\documentclass[12pt]{article} 
\usepackage{amsmath}
\usepackage{amssymb} 
\usepackage{amscd}
\usepackage{float}
\usepackage{theorem} 
\usepackage{graphics}
\usepackage{graphicx}
\usepackage{pb-diagram}

\def\oplusinf{\mathop{\oplus}} 
\def\otimesinf{\mathop{\otimes}}
\def\beq{\begin{equation}}
\def\eeq{\end{equation}}

\def\im{{\mbox{Im}}}
\def\dim{{\mbox{dim}}}
\def\ker{{\mbox{Ker}}}

\def\cala{{\cal A}}

\def\calm{{\cal M}}

\def\cals{{\cal S}}
\def\calc{{\cal C}} 
\def\cale{{\cal E}}

\def\caly{{\cal Y}}
\def\fraca{{\mathfrak A}}

  \def\fracd{{\mathfrak d}} 
\def\bbbone{\mbox{\rm 1\hspace {-.6em} l}}

\def\Yg{{\mathbf Y}}
\def\Tg{{\mathbf T}}
\def\vg{{\mathbf v}}

\newtheorem{theorem}{THEOREM}
\newtheorem{lemma}{LEMMA} 
\newtheorem{proposition}{PROPOSITION}

\begin{document}

\baselineskip=0.7cm

\begin{center} 
 \thispagestyle{empty}
{\large\bf TENSOR FIELDS OF MIXED YOUNG SYMMETRY}
\vspace{0,3cm}
{\large\bf TYPE AND N-COMPLEXES} 

\end{center} 
\vspace{0.75cm}

\begin{center} Michel DUBOIS-VIOLETTE
\footnote{Laboratoire de Physique Th\'eorique, UMR 8627\\ Universit\'e Paris XI,
B\^atiment 210\\ F-91 405 Orsay Cedex, France\\
Michel.Dubois-Violette$@$th.u-psud.fr}  and 
Marc HENNEAUX\footnote{Physique Th\'eorique et Math\'ematique\\ 
Universit\'e Libre de Bruxelles\\
Campus Plaine C.P. 231\\ B-1050 Bruxelles, Belgique
\\
henneaux$@$ulb.ac.be}

\end{center} \vspace{1cm}

\begin{center} \today \end{center}

\vspace {1cm}

\begin{abstract}
 
We construct $N$-complexes of non completely antisymmetric
irreducible tensor fields on $\mathbb R^D$ which generalize
the usual complex $(N=2)$ of differential forms. Although, for
$N\geq 3$, the generalized cohomology of these $N$-complexes is
non trivial, we prove a generalization of the Poincar\'e lemma.
To that end we use a technique reminiscent of the Green ansatz for 
parastatistics.
Several results which appeared in various contexts are shown
to be particular cases of this generalized Poincar\'e lemma. We
furthermore identify the nontrivial part of the generalized
cohomology. Many of the results presented here were announced in 
\cite{D-VH}.

\end{abstract}

  \vspace{1,5cm}
\noindent LPT-ORSAY 01-10\\
ULB-TH/01-16

\newpage

\section{Introduction}

Our aim in this paper is to develop differential tools for
irreducible tensor fields on $\mathbb R^D$ which generalize the
calculus of differential forms. By an irreducible tensor field on
$\mathbb R^D$, we here mean, a smooth mapping $x\mapsto T(x)$ of
$\mathbb R^D$ into a vector space of (covariant) tensors of given
Young symmetry. We recall that this implies that the
representation of $GL_D$ in the corresponding space of tensors is
irreducible.\\

Throughout the following $(x^\mu)=(x^1,\dots,x^D)$ denotes the
canonical coordinates of $\mathbb R^D$ and $\partial_\mu$ are the
corresponding partial derivatives which we identify with the
corresponding covariant derivatives associated to the canonical
flat torsion-free linear connection $\stackrel{(0)}{\nabla}$ of $\mathbb R^D$. 
Thus, for instance, if $T$
is a covariant tensor field of degree $p$ on $\mathbb R^D$ with
components $T_{\mu_1\dots\mu_p}(x)$, then $\stackrel{(0)}{\nabla} T$ denotes
the covariant tensor field of degree $p+1$ with components
$\partial_{\mu_{p+1}}T_{\mu_1\dots\mu_p}(x)$. The operator
$\stackrel{(0)}{\nabla}$ is a first-order differential operator which increases
by one the tensorial degree.\\

In this context, the space $\Omega(\mathbb R^D)$ of differential
forms on $\mathbb R^D$ is the graded vector space of (covariant)
antisymmetric tensor fields on $\mathbb R^D$ with graduation
induced by the tensorial degree whereas the exterior differential
$d$ is up to a sign the composition of the above $\stackrel{(0)}{\nabla}$ with
antisymmetrisation, i.e.
\begin{equation}
d=(-1)^p{\mathbf A}_{p+1}\circ \stackrel{(0)}{\nabla} : \Omega^p(\mathbb R^D)\rightarrow
\Omega^{p+1}(\mathbb R^D)
\label{eq1}
\end{equation}
where ${\mathbf A}_p$ denotes the antisymmetrizer on tensors of degree 
$p$. The sign factor $(-1)^p$ arises because $d$ acts from the left, 
while we defined $(\stackrel{0}{\nabla}T)_{\mu_{1}\dots \mu_{p+1}}= 
\partial_{\mu_{p+1}}T_{\mu_{1}\dots\mu_{p}}$. One has  $d^2=0$ and
the Poincar\'e lemma asserts that the cohomology of the complex
$(\Omega(\mathbb R^D),d)$ is trivial, i.e. that one has
$H^p(\Omega(\mathbb R^D))=0$,  $\forall p\geq 1$ and $H^0(\Omega(\mathbb 
R^D))= \mathbb R$ where $H(\Omega(\mathbb 
R^D))=\ker(d)/\im(d)=\oplusinf_{p}H^p(\Omega(\mathbb R^D))$ 
with $H^p(\Omega(\mathbb R^D))=\ker(d:\Omega^p(\mathbb R^D)\rightarrow
\Omega^{p+1}(\mathbb R^D))/d(\Omega^{p-1}(\mathbb R^D))$.\\

From the point of view of Young symmetry, antisymmetric
tensors correspond to Young diagrams (partitions) described 
by one column of cells, corresponding to the partition $(1^p)$,
whereas ${\mathbf A}_p$ is the associated Young symmetrizer,
(see next section for definitions and conventions).\\

There is a relatively easy way to generalize the pair
$(\Omega(\mathbb R^D),d)$ which we now describe. Let
$(Y)=(Y_p)_{p\in \mathbb N}$ be a sequence of Young diagrams
such that the number of cells of $Y_p$ is $p$, $\forall p\in
\mathbb N$ (i.e. such that $Y_p$ is a partition of the integer
$p$ for any $p$). We define $\Omega^p_{(Y)}(\mathbb R^D)$ to be
the vector space of smooth covariant tensor fields of degree
$p$ on $\mathbb R^D$ which have the Young symmetry type $Y_p$
and we let $\Omega_{(Y)}(\mathbb R^D)$ be the graded vector space
$\displaystyle{\oplusinf_p}\Omega^p_{(Y)}(\mathbb R^D)$. We then
generalize the exterior differential by setting 
\begin{equation}
d=(-1)^p{\mathbf Y}_{p+1}\circ\stackrel{(0)}{\nabla} :\Omega^p_{(Y)}(\mathbb
R^D)\rightarrow \Omega^{p+1}_{(Y)}(\mathbb R^D)
\end{equation}
where ${\mathbf Y}_p$ is now the Young symmetrizer on tensor
of degree $p$ associated to the Young symmetry $Y_p$. This $d$
is again a first order differential operator which is of
degree one, (i.e. it increases the tensorial degree by one),
but now, $d^2\not= 0$ in general. Instead, one has the
following result.

\begin{lemma}
Let $N$ be an integer with $N\geq 2$ and assume that $(Y)$ is
such that the number of columns of the Young diagram $Y_p$ is
strictly smaller than $N$ (i.e. $\leq N-1$) for any $p\in
\mathbb N$. Then one has $d^N=0$.
\end{lemma}

In fact the indices in one column are antisymmetrized (see
below) and $d^N\omega$ involves necessarily at least two
partial derivatives $\partial$ in one of the columns since
there are $N$ partial derivatives involved and at most $N-1$
columns.\\

Thus if $(Y)$ satisfies the condition of Lemma 1, the pair 
$(\Omega_{(Y)}(\mathbb R^D),d)$ is a $N$-{\sl complex} (of cochains)
\cite{Kap}, \cite{D-V}, \cite{D-VK}, \cite{KW}, \cite{D-V2},
i.e. here a graded vector space equipped with an endomorphism
$d$ of degree 1, its $N$-{\sl differential}, satisfying $d^N=0$. 
Concerning $N$-complexes, we
shall use here the notations and the results  of \cite{D-V2} which
will be recalled when needed.\\

Notice that $\Omega^p_{(Y)}(\mathbb R^D)=0$ if the first column of
$Y_p$ contains more than $D$ cells and that therefore, if $Y$
satisfies the condition of Lemma 1, then $\Omega^p_{(Y)}(\mathbb
R^D)=0$ for $p>(N-1)D$. \\

One can also define a graded bilinear product on
$\Omega_{(Y)}(\mathbb R^D)$ by setting 
\begin{equation}
(\alpha\beta)(x)={\mathbf Y}_{a+b}(\alpha(x)\otimes \beta(x))
\label{eq3}
\end{equation}
for $\alpha\in \Omega^a_{(Y)}(\mathbb R^D)$, $\beta\in
\Omega^b_{(Y)}(\mathbb R^D)$ and $x\in \mathbb R^D$. This product is
by construction bilinear with respect to the $C^\infty(\mathbb
R^D)$-module structure of $\Omega_{(Y)}(\mathbb R^D)$ (i.e. with
respect to multiplication by smooth functions). It is worth
noticing here that one always has $\Omega^0_{(Y)}(\mathbb
R^D)=C^\infty(\mathbb R^D)$.

In this paper we shall not stay at this level of generality;
for each $N\geq 2$ we shall choose a maximal $(Y)$, denoted by
$(Y^N)=(Y^N_p)_{p\in\mathbb N}$, satisfying the condition of lemma
1. The Young diagram with $p$ cells $Y^N_p$ is defined in the
following manner: write the division of $p$ by $N-1$, i.e. write
$p=(N-1)n_p+r_p$ where $n_p$ and $r_p$ are (the unique) integers
with $0\leq n_p$ and $0\leq r_p\leq N-2$ ($n_p$ is the quotient
whereas $r_p$ is the remainder), and let $Y^N_p$ be the Young
diagram with $n_p$ rows of $N-1$ cells and the last row with
$r_p$ cells (if $r_p\not= 0$). One has
$Y^N_p=((N-1)^{n_p},r_p)$, that is we fill the rows maximally.

\begin{figure}[H]
\begin{center}
\includegraphics{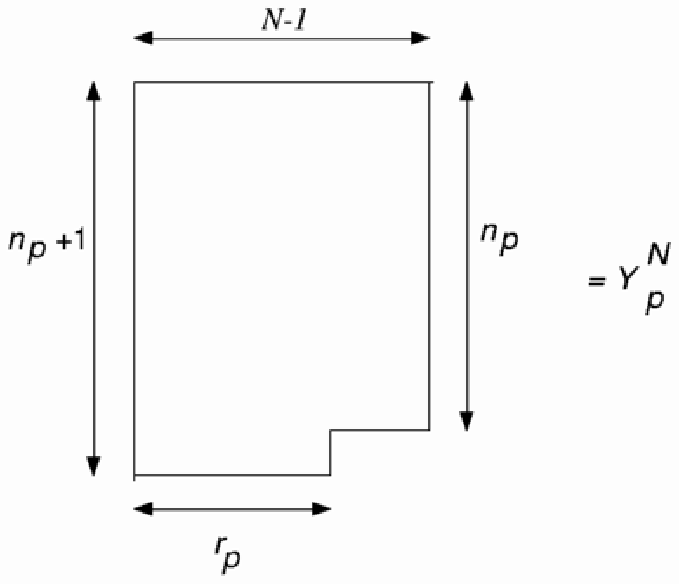}
\end{center}
\end{figure}

\noindent We shall denote $\Omega_{(Y^N)}(\mathbb R^D)$ and
$\Omega^p_{(Y^N)}(\mathbb R^D)$ by $\Omega_N(\mathbb R^D)$ and
$\Omega^p_N(\mathbb R^D)$, respectively. It is clear that $(\Omega_2(\mathbb
R^D),d)$ is the usual complex $(\Omega(\mathbb R^D),d)$ of
differential forms on $\mathbb R^D$. The $N$-complex
$(\Omega_N(\mathbb R^D),d)$ will be simply denoted by
$\Omega_N(\mathbb R^D)$. We recall \cite{D-V2} that the
(generalized) cohomology of the $N$-complex $\Omega_N(\mathbb
R^D)$ is the family of graded vector spaces $H_{(k)}(\Omega_N(\mathbb 
R^D))$ $k\in
\{1,\dots,N-1\}$  defined by
$H_{(k)}(\Omega_N(\mathbb R^D))=\ker(d^k)/\im(d^{N-k})$, i.e.
$H_{(k)}(\Omega_N(\mathbb
R^D))=\displaystyle{\oplusinf_p}H^p_{(k)}(\Omega_N(\mathbb
R^D))$ with 
\[
H^p_{(k)}(\Omega_N(\mathbb
R^D))=\ker(d^k:\Omega^p_N(\mathbb R^D)\rightarrow \Omega^{p+k}_N
(\mathbb R^D))/d^{N-k}(\Omega^{p+k-N}(\mathbb R^D)).
\]
The following statement is our generalization of the Poincar\'e
lemma.

\begin{theorem}\label{theo1}
One has $H^{(N-1)n}_{(k)}(\Omega_N(\mathbb R^D))=0$, $\forall
n\geq 1$ and $H^0_{(k)}(\Omega_N(\mathbb R^D))$ is the space of
real polynomial functions on $\mathbb R^D$ of degree strictly
less than $k$ (i.e. $\leq k-1$) for $k\in\{1,\dots,N-1\}$.
\end{theorem}

This statement reduces to the Poincar\'e lemma for $N=2$ but it
is a nontrivial generalization for $N\geq 3$ in the sense that,
as we shall see, the spaces $H^p_{(k)}(\Omega_N(\mathbb R^D))$
are nontrivial for $p\not=(N-1)n$ and, in fact, are generically
infinite dimensional for $D\geq 3$, $p\geq N$.\\

The connection between the complex of differential forms on
$\mathbb R^D$ and the theory of classical gauge field of spin 1
is well known. Namely the subcomplex
\begin{equation}
\Omega^0(\mathbb R^D)\stackrel{d}{\rightarrow}\Omega^1(\mathbb
R^D)\stackrel{d}{\rightarrow}\Omega^2(\mathbb R^D)
\stackrel{d}{\rightarrow}\Omega^3(\mathbb R^D)
\label{eq4}
\end{equation}
has the following interpretation in terms of spin 1 gauge field
theory. The space $\Omega^0(\mathbb R^D)(=C^\infty(\mathbb R^D))$
is the space of infinitesimal gauge transformations, the space
$\Omega^1(\mathbb R^D)$ is 
the space of gauge potentials (which are the appropriate
description of spin 1 gauge fields to introduce local
interactions). The subspace $d\Omega^0(\mathbb R^D)$ of
$\Omega^1(\mathbb R^D)$ is the space of pure gauge configurations
(which are physically irrelevant), $d\Omega^1(\mathbb R^D)$ is
the space of field strengths or curvatures of gauge potentials.
The identity  $d^2=0$ ensures that the curvatures do not see the
irrelevant pure gauge potentials whereas, at this level, the
Poincar\'e lemma ensures that it is only these irrelevant
configurations which are forgotten when one passes from gauge
potentials to curvatures (by applying $d$). Finally $d^2=0$ also
ensures that curvatures of gauge potentials satisfy the Bianchi
identity, i.e. are in $\ker(d:\Omega^2(\mathbb R^D)\rightarrow
\Omega^3(\mathbb R^D))$, whereas at this level the Poincar\'e
lemma implies that conversely the Bianchi identity characterizes
the elements of $\Omega^2(\mathbb R^D)$ which are curvatures of
gauge potentials.\\

Classical spin 2 gauge field theory is the linearization of
Einstein geometric theory. In this case, the analog of
(\ref{eq4}) is a complex
$\cale^1\stackrel{d_1}{\rightarrow}\cale^2\stackrel{d_2}{\rightarrow}
\cale^3\stackrel{d_3}{\rightarrow}\cale^4$ where $\cale^1$ is the
space of covariant vector field $(x\mapsto X_\mu(x))$ on $\mathbb
R^D$, $\cale^2$ is the space of covariant symmetric tensor fields
of degree 2 ($x\mapsto h_{\mu\nu}(x))$ on $\mathbb R^D$,
$\cale^3$ is the space of covariant tensor fields of degree 4
$(x\mapsto R_{\lambda\mu,\rho\nu}(x))$ on $\mathbb R^D$ having the 
symmetries of
the Riemann curvature tensor and where $\cale^4$
is the space of covariant tensor fields of degree 5 on $\mathbb
R^D$ having the symmetries of the left-hand side of the Bianchi
identity. The arrows $d_1, d_2, d_3$ are given by
\[
\begin{array}{l}
(d_1X)_{\mu\nu}(x)=\partial_\mu X_\nu(x)+\partial_\nu X_\mu (x)\\
\\
(d_2h)_{\lambda\mu,\rho\nu}(x)=\partial_\lambda\partial_\rho 
h_{\mu\nu}(x)
+\partial_\mu\partial_\nu h_{\lambda\rho}(x)-\partial_\mu
\partial_\rho h_{\lambda\nu}(x) -\partial_\lambda\partial_\nu
h_{\mu\rho}(x)\\
\\
(d_3R)_{\lambda\mu\nu,\alpha\beta}(x)=\partial_\lambda
R_{\mu\nu,\alpha\beta}(x)+\partial_\mu
R_{\nu\lambda,\alpha\beta}(x)+\partial_\nu
R_{\lambda\mu,\alpha\beta}(x).
\end{array}
\]
The symmetry of $x\mapsto R_{\lambda\mu,\rho\nu}(x)$,
$\left(\ \begin{tabular}{|c|c|}
\hline
$\lambda$ & $\rho$ \\
\hline
$\mu$ & $\nu$ \\
\hline
\end{tabular}\ 
\right)$, shows that $\cale^3=\Omega^4_3(\mathbb R^D)$ and that
$\cale^4=\Omega^5_3(\mathbb R^D)$; furthermore one canonically
has $\cale^1=\Omega^1_3(\mathbb R^D)$ and
$\cale^2=\Omega^2_3(\mathbb R^D)$. One also sees that $d_1$ and
$d_3$ are proportional to the 3-differential $d$ of
$\Omega_3(\mathbb R^D)$, i.e. $d_1\sim d:\Omega^1_3(\mathbb
R^D)\rightarrow \Omega^2_3(\mathbb R^D)$ and $d_3\sim d :
\Omega^4_3(\mathbb R^D)\rightarrow \Omega^5_3(\mathbb R^D)$.
The structure of $d_2$ looks different, it is of second order
and increases by 2 the tensorial degree. However it is easy to
see that it is proportional to $d^2:\Omega^2_3(\mathbb R^D)
\rightarrow \Omega^4_3(\mathbb R^D)$. Thus the analog of
(\ref{eq4}) is (for spin 2 gauge field theory)
\begin{equation}
\Omega^1_3(\mathbb R^D)\stackrel{d}{\rightarrow}
\Omega^2_3(\mathbb R^D)\stackrel{d^2}{\rightarrow}
\Omega^4_3(\mathbb R^D)\stackrel{d}{\rightarrow}
\Omega^5_3(\mathbb R^D)
\label{eq5}
\end{equation}
and the fact that it is a complex follows from $d^3=0$ whereas
our generalized Poincar\'e lemma (Theorem 1) implies that it is
in fact an exact sequence. Exactness at $\Omega^2_3(\mathbb
R^D)$ is $H^2_{(2)}(\Omega_3(\mathbb R^D))=0$ and exactness
at $\Omega^4_3(\mathbb R^D)$ is $H^4_{(1)}(\Omega_3(\mathbb
R^D))=0$, (the exactness at $\Omega^4_3(\mathbb R^D)$ is the
main statement of \cite{Gas}).\\

Thus what plays the role of the complex of differential forms
for the spin 1 (i.e. $\Omega_2(\mathbb R^D))$ is the 3-complex
$\Omega_3(\mathbb R^D)$ for the spin 2. More generally, for the
spin $S\in \mathbb N$, this role is played by the
$(S+1)$-complex $\Omega_{S+1}(\mathbb R^D)$. In particular, the
analog of the sequence (\ref{eq4}) for the spin 1 is the
complex
\begin{equation}
\Omega^{S-1}_{S+1}(\mathbb
R^D)\stackrel{d}{\rightarrow}\Omega^S_{S+1}(\mathbb
R^D)\stackrel{d^S}{\rightarrow}\Omega^{2S}_{S+1}(\mathbb
R^D)\stackrel{d}{\rightarrow} \Omega^{2S+1}_{S+1}(\mathbb R^D)
\label{eq6}
\end{equation}
for the spin $S$. The fact that (\ref{eq6}) is a complex was
known, \cite{dWF}, it here follows from $d^{S+1}=0$. One easily
recognizes that $d^S:\Omega^S_{S+1}(\mathbb R^D)\rightarrow
\Omega^{2S}_{S+1}(\mathbb R^D)$ is the generalized (linearized)
curvature of \cite{dWF}. Our theorem 1 implies that sequence
(\ref{eq6}) is exact: exactness at $\Omega^S_{S+1}(\mathbb
R^D)$ is $H^S_{(S)}(\Omega_{S+1}(\mathbb R^D))=0$ whereas
exactness at $\Omega^{2S}_{S+1}(\mathbb R^D)$ is
$H^{2S}_{(1)}(\Omega_{S+1}(\mathbb R^D)=0$, (exactness at 
$\Omega^S_{S+1}(\mathbb R^D)$ was directly proved in \cite{DD} for the 
case $S=3$).\\

Finally, there is a generalization of Poincar\'e duality for
$\Omega_N(\mathbb R^D)$, which is obtained by contractions of
the columns with the Kroneker tensor
$\varepsilon^{\mu_1\dots\mu_D}$ of $\mathbb R^D$, that we shall
describe in this paper. When combined with Theorem 1, this
duality leads to another kind of results. A typical result of
this kind is the following one. Let $T^{\mu\nu}$ be a symmetric
contravariant tensor field of degree 2 on $\mathbb R^D$
satisfying $\partial_\mu T^{\mu\nu}=0$, (like e.g. the stress
energy tensor), then there is a contravariant tensor field
$R^{\lambda\mu\rho\nu}$ of degree 4 with the symmetry
$ \begin{tabular}{|c|c|}
\hline
$\lambda$ & $\rho$ \\
\hline
$\mu$ & $\nu$ \\
\hline
\end{tabular}\ 
$, (i.e. the symmetry of Riemann curvature tensor), such
that
\begin{equation}
T^{\mu\nu}=\partial_\lambda\partial_\rho R^{\lambda\mu \rho\nu}
\label{eq7}
\end{equation}
In order to connect this result with Theorem 1, define
$\tau_{\mu_1\dots\mu_{D-1}
\nu_1\dots\nu_{D-1}}=\linebreak[4]T^{\mu\nu}\varepsilon_{\mu\mu_1\dots\mu_{D-1}}
\varepsilon_{\nu\nu_1\dots \nu_{D-1}}$. Then one has
$\tau\in\Omega^{2(D-1)}_3(\mathbb R^D)$ and conversely, any
$\tau\in \Omega^{2(D-1)}_3(\mathbb R^D)$ can be expressed in
this form in terms of a symmetric contravariant 2-tensor. It is
easy to verify that $d\tau=0$ (in $\Omega_3(\mathbb R^D))$ is
equivalent to $\partial_\mu T^{\mu\nu}=0$. On the other hand,
Theorem 1 implies that $H^{2(D-1)}_{(1)}(\Omega_3(\mathbb
R^D))=0$ and therefore $\partial_\mu T^{\mu\nu}=0$ implies that
there is a $\rho\in \Omega^{2(D-2)}_3(\mathbb R^D)$ such that
$\tau=d^2\rho$. The latter is equivalent to (\ref{eq7}) with
$R^{\mu_1\mu_2\ \nu_1\nu_2}$ proportional to
$\varepsilon^{\mu_1\mu_2\dots\mu_D}
\varepsilon^{\nu_1\nu_2\dots
\nu_D}\rho_{\mu_3\dots\mu_D\nu_3\dots\nu_D}$ and one verifies
that, so defined, $R$ has the correct symmetry. That symmetric
tensor fields identically
fulfilling $\partial_\mu T^{\mu\nu}=0$ can
be rewritten as in Eq. (\ref{eq7}) has been used in
\cite{W} and more recently in \cite{BDGH} in the investigation of the consistent deformations of
the free spin two gauge field action.\\

Beside their usefulness for computations (and for unifying various 
results) through the generalization of Poincar\'e lemma (Theorem 
1) and the generalization of the Poincar\'e duality, the 
$N$-complexes described in this paper give a class of nontrivial 
examples of $N$-complexes which are not related with simplicial 
modules. Indeed most nontrivial examples of $N$-complexes considered 
in \cite{D-V}, \cite{D-V2}, \cite{D-V3}, \cite{Kap}, \cite{May}, 
\cite{KW} are of simplicial type and it was shown in \cite{D-V2} that 
such $N$-complexes compute the ordinary (co)homologies of the 
simplicial modules (see also in \cite{KW} for the Hochschild case). 
Furthermore that kind of results have been recently extended to the 
cyclic context in \cite{Wam} where new proofs of above results have 
been carried over. This does not mean that $N$-complexes associated 
with simplicial modules are not useful; for instance in \cite{D-VT2} 
such a $N$-complex (related with a simplicial Hochschild module) was 
needed for the construction of a natural generalized BRS-theory 
\cite{BRS}, \cite{HT} for the zero modes of the $SU(2)$ WZNW-model, 
see in \cite{D-V4} for a general review. It is however very desirable 
to produce useful examples which are not of simplicial type and, apart 
from the universal construction of \cite{D-VK} (and some 
finite-dimensional examples \cite{D-V2}, \cite{D-VK}), the examples produced 
here are the first ones escaping from the simplicial frame.\\

Many results of this paper where announced in our letter \cite{D-VH} so an 
important part of it is devoted to the proofs of these results in 
particular to the proof of Theorem 1 above which generalizes the 
Poincar\'e lemma. In order that the paper be self contained we recall 
some basic definitions and results on Young diagrams and representations of the 
linear group which are needed here. Throughout the paper, we work in 
the real setting, so all vector spaces are on the field $\mathbb R$ of 
real numbers (this obviously generalizes to any commutative field 
$\mathbb K$ of characteristic zero).\\

The plan of the paper is the following. After this introduction we 
discuss Young diagrams, Young symmetry types for tensor and we define 
in this context a notion of contraction. Section 3 is devoted to the 
construction of the basic $N$-complex of tensor fields on $\mathbb 
R^D$ considered in this paper, namely $\Omega_{N}(\mathbb R^D)$, and 
the description of the generalized Poincar\'e (Hodge) duality in this 
context. In Section 4 we introduce a multicomplex on $\mathbb R^D$ 
and we analyse its cohomological properties; Theorem 2 proved there, 
which is by itself of interest, will be the basic ingredient in the 
proof of our generalization of the Poincar\'e lemma i.e. of Theorem 1. 
Section 5 contains this proof of Theorem 1. In Section 6 we analyse 
the structure of the generalized cohomology of $\Omega_{N}(\mathbb 
R^D)$ in the degrees which are not exhausted by Theorem 1. The 
$N$-complex $\Omega_{N}(\mathbb R^D)$ is a generalization of the 
complex $\Omega(\mathbb R^D)=\Omega_{2}(\mathbb R^D)$ of differential 
forms on $\mathbb R^D$; in Section 7 we define another 
generalization $\Omega_{[N]}(\mathbb R^D)$ of the complex of 
differential forms which is also a $N$-complex and which is an 
associative graded algebra acting on the graded space 
$\Omega_{N}(\mathbb R^D)$. In Section 8 which plays the role of a conclusion 
we sketch another possible proof of Theorem 1 based on a 
generalization of algebraic homotopy for $N$-complexes. In this 
section we also define natural $N$-complexes of tensor fields on 
complex manifolds which generalize the usual $\bar\partial$-complex 
(of forms in $d\bar z$).

\section{Young diagrams and tensors}

 For the Young diagrams etc. we use throughout the conventions of 
 \cite{Ful}. A {\sl Young diagram} $Y$ is a diagram which consists of 
 a finite number $r>0$ of rows of identical squares (refered to as the 
 {\sl cells})
 of finite decreasing lengths $m_{1}\geq m_{2}\geq \dots \geq 
 m_{r}>0$ which are arranged with their left hands under one another. 
 The lengths $\tilde m_{1},\dots, \tilde m_{c}$ of the columns of 
 $Y$ are also decreasing $\tilde m_{1}\geq \dots \geq \tilde 
 m_{c}>0$ and are therefore the rows of another Young diagram $\tilde 
 Y$ with $\tilde r=c$ rows. The Young diagram $\tilde Y$ is obtained 
 by flipping $Y$ over its diagonal (from upper left to lower right) 
 and is refered to as the {\sl conjugate} of $Y$. Notice that one 
 has $\tilde m_{1}=r$ and therefore also $m_{1}=\tilde r=c$ and that 
 $m_{1}+\dots +m_{r}=\tilde m_{1}+\dots +\tilde m_{c}$ is the total 
 number of cells of $Y$ which will be denoted by $\vert Y\vert$. It is 
 convenient to add the empty Young diagram $Y_{0}$ characterized by 
 $\vert Y\vert=0$.The figure below describes a Young diagram $Y$ and 
 its conjugate $\tilde Y$ :
 \vspace{0,75cm}
 
 \begin{center}
   $Y$ =\hspace{0,2cm} \begin{tabular}{|c|c|c|c|}\hline
  &  &  &  \\ \hline
   &  &  &  \\ \hline
   &  &  &  \\ \hline
  &  &  &   \\ \hline
    &  & \\ 
 \cline{1-3}
  & \\ 
 \cline{1-2}
  &  \\ 
 \cline{1-2}
 \\
 \cline{1-1}
 \end{tabular}
  \hspace{1cm} $\tilde Y$= \hspace{0,2cm}
 \begin{tabular}{|c|c|c|c|c|c|c|c|}\hline
  & & & & & & & \\ \hline
  & & & & & & \\
  \cline{1-7}
   & & & & \\
   \cline{1-5}
   & & & \\
    \cline{1-4}
    \end{tabular}
\end{center}
  \vspace{1cm}
  
 In the following $E$ denotes a finite-dimensional vector space of 
 dimension $D$ and $E^\ast$ denotes its dual. The $n$-th tensor power 
 $E^{\otimes^n}$ of $E$ identifies canonically with the space of 
 multilinear forms on $(E^\ast)^n$. Let $Y$ be a Young diagram and 
 let us consider that the $\vert Y\vert$ copies of $E^\ast$ in 
 $(E^\ast)^{\vert Y\vert}$ are labelled by the cells of $Y$ so that 
 an  element of $(E^\ast)^{\vert Y\vert}$ is given by specifying an 
 element of $E^\ast$ for each cell of $Y$. The {\sl Schur module} 
 $E^Y$ is defined to be the vector space of all multilinear forms 
 $T$ on $(E^\ast)^{\vert Y\vert}$ such that:
\begin{quote}
 
$(i)$ $T$ is completely antisymmetric in the entries of each 
column of $Y$,

$(ii)$ complete antisymmetrization of $T$ in the entries of a 
column of $Y$ and another entry of $Y$ which is on the right-hand 
side of the column vanishes.
\end{quote}

Notice that $E^Y=0$ if the first column 
of $Y$ has length $\tilde m_{1}>D$. One has $E^Y\subset 
E^{\otimes^{\vert Y\vert}}$ and $E^Y$ is an invariant subspace for 
the action of $GL(E)$ on $E^{\otimes^{\vert Y\vert}}$ which is 
irreducible. Furthermore each irreducible subspace of $E^{\otimes^n}$ 
for the action of $GL(E)$ is isomorphic to $E^Y$ with the above action 
of $GL(E)$ for some Young diagram $Y$ with $\vert Y\vert=n$.\\

Let $Y$ be a Young diagram and let $T$ be an arbitrary 
multilinear form on $(E^\ast)^{\vert Y\vert}$ , ($T\in 
E^{\otimes^{\vert Y\vert}}$). Define the multilinear form 
$\caly(T)$ on $(E^\ast)^{\vert Y \vert}$ by 
\[
\caly(T)=\sum_{p\in R} \sum_{q\in C} 
(-1)^{\varepsilon(q)}T \circ p \circ q
\]
where $C$ is the group of the permutations which permute the
entries of each column and $R$ is the group of the permutations which 
permute the entries of each  row of $Y$. One has $\caly(T)\in 
E^Y$ and the endomorphism $\caly$ of $E^{\otimes^{\vert Y\vert}}$ 
satisfies $\caly^2=\lambda\caly$ for some number $\lambda\not= 0$. 
Thus $\Yg = \lambda^{-1}\caly$ is a projection of $E^{\otimes^{\vert 
Y\vert}}$ into itself, $\Yg^2=\Yg$, with image $\im(\Yg)=E^Y$. The 
projection $\Yg$ will be refered to as the {\sl Young symmetrizer} 
(relative to $E$) of the Young diagram $Y$. The element 
$e_{Y}=\lambda^{-1}\sum_{p\in R}\sum_{q\in C}(-1)^{\varepsilon(q)}pq$ 
of the group algebra of the group $\cals_{\vert Y\vert}$ of 
permutation of $\{1,\dots,\vert Y\vert \}$ is an idempotent which will 
be refered to as the {\sl Young idempotent} of $Y$.\\

By composition of $\Yg$ as 
above with the canonical multilinear mapping of $E^{\vert Y\vert}$ into 
$E^{\otimes^{\vert Y\vert}}$ one obtains a multilinear mapping $\vg 
\mapsto \vg^Y$ of $E^{\vert Y\vert}$ into $E^Y$. The Schur module $E^Y$ together 
with the mapping $\vg \mapsto \vg^Y$ are characterized uniquely up to 
an isomorphism by the following universal property: For any 
multilinear mapping $\phi:E^{\vert Y\vert}\rightarrow F$ of $E^{\vert 
Y\vert}$ into a vector space $F$ satisfying
\begin{quote}
 $(i)$ $\phi$ is completely antisymmetric in the entries of each 
 column of $Y$,
 
 $(ii)$ complete antisymmetrization of $\phi$ in the entries of a 
 column of $Y$ and another entry of $Y$ which is on the right-hand 
 side of the column vanishes,
 \end{quote}
 there is a unique linear mapping 
 $\phi^Y:E^Y\rightarrow F$ such that $\phi(\vg)=\phi^Y(\vg^Y)$. By 
 construction $\vg\mapsto \vg^Y$ satisfies the conditions $(i)$ and 
 $(ii)$ above.\\
 
 There is an obvious notion of inclusion for Young diagrams namely 
 $Y'$ is included in $Y$, $Y'\subset Y$, if one has this inclusion for 
 the corresponding subsets of the plane whenever their upper left 
 cells coincide. This means for instance that $Y'\subset Y$ whenever 
 the length $c=m_{1}$ of the first row of $Y$ is greater than the length 
 $c'=m'_{1}$ of the first row of $Y'$ and that for any $1\leq i\leq 
 c'$ the length $\tilde m_{i}$ of the $i$-th column of $Y$ is greater 
 than the length $\tilde m'_{i}$ of the $i$-th column of $Y'$, ($c\geq 
 c'$ and $\tilde m_{i}\geq \tilde m'_{i}$ for $1\leq i\leq c'$).\\
 
 In the following we shall need a stronger notion. A Young diagram 
 $Y'$ is {\sl strongly included} in another one $Y$ and we write 
 $Y'\subset\subset Y$ if the length of the first row of $Y$ is 
 greater than the length of the first row of $Y'$ and if the length 
 of the {\sl last} column of $Y$ is greater than the length of the 
 first column of $Y'$. Notice that this relation is not reflexive, one 
 has $Y\subset\subset Y$ if and only if $Y$ is {\sl rectangular} which 
 means that all its columns have the same length or equivalently all 
 its rows have the same length. It is clear that $Y'\subset\subset Y$ 
 implies $Y'\subset Y$. \\
 
 Let $Y$ and $Y'$ be Young diagrams such that $Y'\subset \subset Y$ 
 and let $\tilde m_{1}\geq \dots \geq \tilde m_{c}>0$ be the lengths 
 of the columns of $Y$ and $\tilde m'_{1}\geq \dots \geq \tilde 
 m'_{c'}>0$ be the lengths of the columns of $Y'$; one has $c\geq c'$ 
 and $\tilde m_{c}\geq \tilde m'_{1}$. Define {\sl the contraction 
 of $Y$ by $Y'$} to be the Young diagram $\calc(Y\vert Y')$ obtained 
 from $Y$ by dropping $\tilde m'_{1}$ cells of the last i.e. the 
 $c$-th column of $Y$, $\tilde m'_{2}$ cells of the ($c-1$)-th column 
 of $Y,\dots,\tilde m'_{c'}$ cells of the ($c-c'+1$)-th column of 
 $Y$. If $\tilde m_{c}$ is strictly geater than $\tilde m'_{1}$ then 
 $\calc(Y\vert Y')$ has $c$ columns as $Y$, however if $\tilde 
 m_{c}=\tilde m'_{1}$ then the number of columns of $\calc(Y\vert Y')$ 
 is strictly smaller than $c$ (it is $c-1$ if $\tilde m_{c-1}$ is 
 strictly greater than $\tilde m'_{2}$, etc.). Notice that if $Y$ is 
 rectangular then $\calc(Y\vert Y')\subset\subset Y$ and 
 $\calc(Y\vert\calc(Y\vert Y'))=Y'$ so that $Y'\mapsto \calc (Y\vert Y')$ 
 is then an involution on the set of Young diagrams $Y'$ which are 
 strongly included in $Y$ ($Y'\subset \subset Y$).\\
 
 Let again $Y$ and $Y'$ be Young diagrams with $Y'\subset\subset Y$. 
 Our aim is now to define a bilinear mapping $(T,T')\mapsto 
 \calc(T\vert T')$ of $E^Y\times E^{\ast Y'}$ into $E^{\calc(Y\vert Y')}$. 
 This will be obtained by restriction of a bilinear mapping 
 $(T,T')\mapsto\calc(T\vert T')$ of $E^{\otimes^{\vert Y\vert}}\times 
 E^{\ast \otimes^{\vert Y'\vert}}$ into $E^{\otimes^{\vert\calc(Y\vert 
 Y')\vert}}$ which will be an ordinary (complete) tensorial 
 contraction. Any such tensorial contraction associates to a 
 contravariant tensor $T$ of degree $\vert Y\vert$ (i.e. $T\in 
 E^{\otimes^{\vert Y\vert}}$) and a covariant tensor $T'$ of degree 
 $\vert Y'\vert$ (i.e. $T'\in E^{\ast\otimes^{\vert Y'\vert}}$) a 
 contravariant tensor of degree $\vert \calc(Y\vert Y')\vert$, 
 ($Y'\subset\subset Y$). In order to specify such a contraction, one 
 has to specify the entries of $T$, that is of $Y$, to which each 
 entry of $T'$, that is of $Y'$, is contracted (recalling that $T$ is 
 a linear combination of canonical images of elements of $E^{\vert 
 Y\vert}$ and that $T'$ is a linear combination of canonical images 
 of elements of $E^{\ast\vert Y'\vert}$). In order that $\calc(T\vert 
 T')$ has the right antisymmetry in the entries of each column of 
 $\calc(Y\vert Y')$ when $T\in E^Y$ and $T'\in E^{\ast Y'}$, one has to 
 contract the entries of $T'$ corresponding to the $i$-th column of 
 $Y'$ with entries of $T$ corresponding to the ($c-i+1$)-th column 
 of $Y$. The precise choice and the order of the latter entries is 
 irrelevant up to a sign in view of the antisymmetry in the entries of 
 a column. Our choice is to contract the first entry of the $i$-th 
 column of $Y'$ with the last entry of the ($c-i+1$)-th column of 
 $Y$, the second entry of the $i$-th column of $Y'$ with the 
 penultimate entry of the ($c-i+1$)-th column of $Y$, etc. for any 
 $1\leq i\leq c'$ (with obvious conventions). This fixes the bilinear 
 mapping $(T,T')\mapsto \calc(T\vert T')$ of $E^{\otimes^{\vert 
 Y\vert}}\times E^{\ast\otimes^{\vert Y'\vert}}$ into 
 $E^{\otimes^{\vert \calc(Y\vert Y')\vert}}$. The following figure 
 describes picturally in a particular case the construction of 
 $\calc(Y\vert Y')$ as well as the places where the contractions are 
 carried over in the corresponding construction of $\calc(T\vert 
 T')$~:
 
 \vspace{1cm}
 
  $Y$ =\hspace{0,2cm} \begin{tabular}{|c|c|c|c|}\hline
  &  &  &  \\ \hline
   &  &  &  \\ \hline
   &  &  &  \\ \hline
  &  &  &   \\ \hline
    &  & \\ 
 \cline{1-3}
  & \\ 
 \cline{1-2}
  &  \\ 
 \cline{1-2}
 \\
 \cline{1-1}
 \end{tabular}
  \hspace{0,2cm}
  $\longrightarrow$ \hspace{0,2cm}
 \begin{tabular}{|c|c|c|c|}\hline
  &  &  &  \\ \hline
   &  &  &  \\ \hline
   &  &  &  \\ \hline
  &  &  &   \\ \hline
    &  & \\ 
 \cline{1-3}
  & \\ 
 \cline{1-2}
  &  \\ 
 \cline{1-2}
 \\
 \cline{1-1}
 \end{tabular} 
  \hspace{0,2cm}
  $\longrightarrow$ \hspace{0,2cm}
\begin{tabular}{|c|c|c|c|}\hline
  &  &  &  \\ \hline
   &  &  &  \\ \hline
   &  &    \\ \cline{1-3}
  &  &     \\ \cline{1-3}
    &  \\ \cline{1-2}
  & \\ 
 \cline{1-2}
  &  \\ 
 \cline{1-2}
 \\
 \cline{1-1}
 \end{tabular} \hspace{0,2cm} = \hspace{0,2cm} $\calc(Y\vert Y')$\\

 \hspace{1,1cm} \raisebox{-2mm}{ $Y'$ = \begin{tabular}{|c|c|}\hline
 & \\ \hline
 \\
 \cline{1-1}
 \end{tabular}} \hspace{1,6cm}
 \rotatebox{180}{\begin{tabular}{|c|c|}\hline
 & \\ \hline
 \\
 \cline{1-1}
  \end{tabular} }\hspace{0,1cm}\raisebox{0,5cm}{$\uparrow$}
  \vspace{1cm}

 \begin{proposition}
  Let $T$ be an element of $E^Y$ and $T'$ be an element of $E^{\ast Y'}$ 
  with $Y'\subset\subset Y$. Then $\calc(T\vert T')$ is an element 
  of $E^{\calc(Y\vert Y')}$.
 \end{proposition}
 
 \noindent {\bf Proof}
As before, we identify $\calc(T\vert T')\in E^{\otimes^{\vert 
\calc(Y\vert Y')|}}$ with a multilinear form on $E^{\ast\vert 
\calc(Y\vert Y')\vert}$. To show that $\calc(T\vert T')$ is in 
$E^{\calc(Y\vert Y')}$ means  verifying properties $(i)$ and $(ii)$ 
above. Property $(i)$, i.e. antisymmetry in the columns entries of 
$\calc(Y\vert Y')$, is clear. Property $(ii)$ has to be verified for 
each column of $\calc(Y\vert Y')$ and entry on its right-hand side 
which can be chosen to be the first entry of a column on the 
right-hand side (in view of the column antisymmetry). If the column 
is the last one it has no entry on the right-hand side so there 
nothing to verify and if the column is a full column of $Y$, i.e. has 
not be contracted, which is the case for the $i$-th column with $i\leq 
c-c'$, the property $(ii)$ follows from the  
same property for $T$ (assumption $T\in E^Y$) . Thus to achieve the proof of the proposition we 
only need to verify property $(ii)$ in the case where both $Y$ and 
$Y'$ have exactly two columns of lengths say $\tilde m_{1}\geq \tilde 
m_{2}$ for $Y$ and $\tilde m'_{1}\geq \tilde m'_{2}$ for $Y'$ with 
$\tilde m_{2}>\tilde m'_{1}$. In this case $\calc(Y\vert Y')$ has 
also two columns of lengths $\tilde m_{1}-\tilde m'_{2}$ and $\tilde 
m_{2}-\tilde m'_{1}$ ($\tilde m_{1}-\tilde m'_{2}\geq \tilde 
m_{2}-\tilde m'_{1}>0$) and one has to verify that antisymmetrization 
of the first entry of the second column of $\calc(Y\vert Y')$  with the entries
of the first column (of length $\tilde m_{1}-\tilde m'_{2}$)
of $\calc(Y\vert Y')$ in $\calc(T\vert T')$ gives zero. We know that 
antisymmetrization with all entries of the first column of $Y$ give 
zero (for $T$); however when contracted with $T'$ this identity 
implies a sum of antisymmetrizations of the entries of the first 
column of $Y'$ with the successive entries of its second column for 
$T'$ which gives zero ($T'=E^{\ast Y'}$) and reduces therefore to 
desired antisymmetrization with the $\tilde m_{1}-\tilde m'_{2}$ 
first entries.$\square$

 \section{Generalized complexes of tensor fields}
 
 Throughout this section $(Y)$ denotes not just one Young diagram but a 
 sequence $(Y)=(Y_{p})_{p\in \mathbb N}$ of Young diagrams $Y_{p}$ such 
 that the number of cells of $Y_{p}$ is equal to $p$ that is $\vert 
 Y_{p}\vert=p$, $\forall p\in \mathbb N$. Notice that there is no 
 freedom for $Y_{0}$ and $Y_{1}$: $Y_{0}$ must be the empty Young 
 diagram and $Y_{1}$ is the Young diagram with one cell. Let us 
 denote by $\wedge_{(Y)}E$ the direct sum $\oplus_{p\in \mathbb N} 
 E^{Y_{p}}$ of the Schur modules $E^{Y_{p}}$. This is a graded vector 
 space with $\wedge^p_{(Y)}E=E^{Y_{p}}$. The  origin of this notation 
 is that for the sequence $(Y^2)=(Y^2_{p})$ of the one column Young 
 diagrams, i.e. $Y^2_{p}$ is the Young diagram with $p$ cells in one 
 column for any $p\in \mathbb N$, then $\wedge_{(Y^2)}E$ is the 
 exterior algebra $\wedge E$ of $E$.\\
    
    In the following, we shall be interested in particular sequences 
    $(Y^N)=(Y^N_{p})_{p\in \mathbb N}$ of Young diagrams satisfying the 
    assumption of Lemma 1 (as explained in the introduction). The 
    sequence $(Y^N)$ contains Young diagrams $Y^N_{p}$ in which all 
    the rows but the last one are of length $N-1$, the last one being 
    of length smaller than or equal to $N-1$ in such a way that $\vert 
    Y^N_{p}\vert=p$ ($\forall p\in \mathbb N$). Picturally one has 
    for instance for $N=5$\\
    
    \vspace{1cm}
    $Y^5_{3}$ = \begin{tabular}{|c|c|c|}\hline
     & & \\\hline
     \end{tabular} \hspace{1cm}
     $Y^5_{22}$ = \begin{tabular}{|c|c|c|c|}\hline
     & & & \\ \hline
     & & & \\ \hline
     & & & \\ \hline
     & & & \\ \hline
     & & & \\ \hline
     & \\ \cline{1-2}
     \end{tabular}\hspace{1cm} $Y^5_{24}$ = \begin{tabular}{|c|c|c|c|}\hline
      & & & \\ \hline
       & & & \\ \hline
        & & & \\ \hline
	 & & & \\ \hline
	  & & & \\ \hline
	   & & & \\ \hline
	   \end{tabular}\\
	   \vspace{1cm}
	   
	   \noindent and so on. In this case $\wedge_{(Y^N)}E$ and 
	   $\wedge^p_{(Y^N)}E=E^{Y^N_{p}}$ will be simply denoted by 
	   $\wedge_{N}E$ and $\wedge^p_{N}E$ respectively. Notice that 
	   $\wedge^p_{N}E=0$ for $p>(N-1)D$, ($D=\dim E$), so that 
	   $\wedge_{N}(E)=\oplus^{(N-1)D}_{p=0}\wedge^p_{N}E$ is 
	   finite-dimensional.\\

	   Let us assume that $E$ is equipped with a dual  
	   volume, i.e. a non-vanishing  element $\varepsilon$ of $\wedge^D E$ 
	   ($=\wedge^D_{2}E$), which is therefore a basis of the 
	   $\mbox{1-dimensional}$ space $\wedge^D (E)$. It is straightforward that 
	   $\varepsilon^{\otimes^{(N-1)}}$ is in 
	   $\wedge^{(N-1)D}_{N}E=E^{Y^N_{(N-1)D}}$ because $(i)$ is obvious 
	   whereas $(ii)$ is trivial i.e. empty. The Young diagram 
	   $Y^N_{(N-1)D}$ is rectangular so that each Young diagram which is 
	   included in $Y^N_{(N-1)D}$ is in fact strongly included in 
	   $Y^N_{(N-1)D}$; this is in particular the case for the $Y^N_{p}$ 
	   for $p\leq (N-1)D$. One then defines a linear isomorphism 
	   $\ast:\wedge_{N}E^\ast\rightarrow \wedge_{N}E$ generalizing the 
	   algebraic part of the Poincar\'e (Hodge) duality by setting
	   \beq
	   \ast \omega = \calc (\varepsilon^{\otimes^{(N-1)}}\vert \omega)
	   \label{eq10}
	   \end{equation}
	   for $\omega\in \wedge_{N}E^\ast$. One has
	   \beq
	   \ast \wedge^p_{N}E^\ast=\wedge^{(N-1)D-p}_{N}E
	   \label{eq11}
	   \end{equation}
	   for $p=0,\dots,(N-1)D$.\\
   
 Let $(e_{\mu})_{\mu\in\{1,\dots,D\}}$ be a basis of $E$ and let 
 $(\theta^\mu)$ be the dual basis of $E^\ast$. Our aim is to be able 
 to compute in terms of the components of tensors for the various 
 concepts connected with Young diagrams. For this, one has to decide 
 the linear order in which one writes the components of a tensor $T\in 
 E^{\otimes^{\vert Y\vert}}$ or, which is the same, of a multilinear 
 form $T$ on $E^{\ast\vert Y\vert}$ for any given Young diagram $Y$. 
 Since we have labelled the arguments (entries) of such a $T$ by the 
 cells of $Y$ and since the components are obtained by taking the 
 arguments among the $\theta^\mu$, this means that one has to choose 
 an order for the cells of $Y$ (i.e a way to ``read the diagram" 
 $Y$). One natural choice is to read the rows of $Y$ from left to 
 right and then from up to down (like a book); another natural choice 
 is to read the columns of $Y$ from up to down and then from left to 
 right. Although the first choice is very natural with respect to the 
 sequences $(Y^N)$ of Young diagrams introduced above and will be used 
 later, we shall choose 
 the second way of ordering in the following. The reason is that when 
 $T$ belongs to the Schur module $E^Y$, then it is (property $(i)$) 
 antisymmetric in the entries of each columns. Thus if $Y$ has columns 
 of lengths $\tilde m_{1}\geq \dots \geq \tilde m_{c}$ ($>0$ for 
 $\vert Y\vert\not= 0$) our choice is induced by the {\sl canonical 
 identification}
 \beq
 E^Y \subset \wedge^{\tilde m_{1}} E \otimes \dots \otimes 
 \wedge^{\tilde m_{c}} E
 \label{Id}
 \end{equation}
 of the Schur module $E^Y$ as a subspace of 
 $\wedge^{\tilde m_{1}} E\otimes \dots \otimes \wedge^{\tilde m_{c}} 
 E$ where $\wedge^p E=\wedge^p_{2}E$ is the $p$-th exterior power of 
 $E$. With the above choice, the components (relative to the basis 
 $(e_{\mu})$ of $E$) of $T\in E^{\otimes^{\vert Y\vert}}$ read 
 $T^{\mu^1_{1}\dots \mu^{\tilde m_{1}}_{1}, \dots ,\mu^1_{c}\dots 
 \mu^{\tilde m_{c}}_c}$ and $T\in E^Y$ if and only if these components 
 are completely antisymmetric in the $\mu^1_{r},\dots,\mu^{\tilde 
 m_{r}}_{r}$ for each $r\in \{1,\dots,c\}$ and such that complete 
 antisymmetrization in the $\mu^1_{r},\dots,\mu^{\tilde m_{r}}_{r}$ 
 and $\mu^1_{s}$ gives zero for any $1\leq r < s\leq c$.\\
 
 We have defined for a sequence $(Y)=(Y_{p})$ of Young diagrams 
 with $\vert Y_{p} \vert =p$ ($\forall p\in \mathbb N$) the  graded 
 vector space
 $\wedge_{(Y)}E$ which can be considered as a generalization of the 
 exterior algebra $\wedge E$ as explained above. We now wish to 
 define the corresponding generalization of differential forms. Let 
 $M$ be a $D$-dimensional smooth manifold. For any Young diagram $Y$ 
 one has the smooth vector bundle $T^{\ast Y} (M)$ over $M$ of the 
 Schur modules $(T^\ast_{x}(M))^Y$, $x\in M$. Correspondingly, for 
 $(Y)$ as above, one has the smooth bundle $\wedge_{(Y)}T^\ast(M)$ 
 over $M$ of  graded vector spaces
 $\wedge_{(Y)}T^\ast_{x}(M)$. 
 The graded $C^\infty(M)$-module $\Omega_{(Y)}(M)$ of smooth sections of 
 $\wedge_{(Y)}T^\ast(M)$ is the 
 generalization  of differential forms 
 corresponding to $(Y)$. In order to generalize the exterior 
 differential one has to choose a connection $\nabla$ on the vector 
 bundle $T^\ast(M)$ that is a linear connection $\nabla$ on $M$. Such 
 a connection extends canonically as linear mappings
 \[
 \nabla:\Omega^p_{(Y)}(M)\rightarrow \Omega^p_{(Y)}(M) 
 \otimesinf_{C^\infty(M)} \Omega^1(M)
 \]
 where $\Omega^1(M)=\Omega^1_{(Y)}(M)$ is the $C^\infty(M)$-module of 
 smooth sections of $T^\ast(M)$ (i.e. of differential 1-forms) 
 satisfying
 \[
 \nabla (\alpha f)=\nabla(\alpha)f + \alpha\otimes df
 \]
 for any $\alpha\in \Omega^p_{(Y)}(M)$ and $f\in C^\infty(M)$ and 
 where $d$ is the ordinary differential of $C^\infty(M)$ into 
 $\Omega^1(M)$. Notice that for any sequence $(Y)$ of Young diagrams 
 as above, one has $\Omega^0_{(Y)}=\Omega^0(M)=C^\infty(M)$ and 
 $\Omega^1_{(Y)}(M)=\Omega^1(M)$ since one has  no choice for $Y_{0}$ 
 and $Y_{1}$. Let us define the generalization of the 
 covariant exterior differential
 $d_{\nabla}:\Omega_{(Y)}(M)\rightarrow \Omega_{(Y)}(M)$  by
 \beq
 d_{\nabla}=(-1)^p \Yg_{p+1}\circ \nabla:\Omega^p_{(Y)}(M)\rightarrow 
 \Omega^{p+1}_{(Y)}(M)
 \label{diff}
 \end{equation}
 for any $p\in \mathbb N$. Notice that $d_{\nabla}=d$ on 
 $C^\infty(M)=\Omega^0_{(Y)}(M)$ and that $d_{\nabla}$ is a first 
 order differential operator. Lemma 1 in the introduction admits the 
 following generalization.
 \begin{lemma}\label{lem1b}
  Let $N$ be an integer with $N\geq 2$ and assume that $(Y)$ is such 
  that the number of columns of the Young diagram $Y_{p}$ is strictly 
  smaller than $N$ for any $p\in \mathbb N$. Then $(d_{\nabla})^N$ is 
  a differential operator of order strictly smaller than $N$. If 
  $\nabla$ is torsion-free, then $d_{\nabla}^N$ is order
  strictly smaller than $N-1$. If 
  furthermore $\nabla$ has vanishing torsion and curvature then one 
  has $(d_{\nabla})^N=0$.
  \end{lemma}
  The proof is straightforward. In the case $N=2$, if $\nabla$ is 
  torsion free, $(d_{\nabla})^2$ is not only an operator of order zero 
  but $(d_{\nabla})^2=0$ follows from the first Bianchi identity; 
  however in this case, for $(Y^2)$, $d_{\nabla}$ coincides with the 
  ordinary exterior differential. For the sequences $(Y^N)=(Y^N_{p})$ 
  we denote $\Omega_{(Y^N)}(M)$ and $\Omega^p_{(Y^N)}(M)$ simply by 
  $\Omega_{N}(M)$ and $\Omega^p_{N}(M)$. As already mentioned 
  $\Omega_{2}(M)$ is the graded algebra $\Omega(M)$ of differential 
  forms on $M$. \\
  
 Not every $M$ admits a flat torsion-free linear connection. 
  In the following we shall concentrate on $\Omega_{N}(\mathbb R^D)$ 
  equipped with $d=d_{\stackrel{(0)}{\nabla}}$ where $\stackrel{(0)}{\nabla}$ is the canonical 
  flat torsion-free connection of $\mathbb R^D$. So equipped, 
  $\Omega_{N}(\mathbb R^D)$ is a $N$-complex. One has of course 
  $\Omega_{N}(\mathbb R^D)=\wedge_{N}\mathbb R^{D\ast}\otimes 
  C^\infty (\mathbb R^D)$. Let us equip $\mathbb R^D$ with the dual 
  volume $\varepsilon \in \wedge^D\mathbb R^D$ which is the 
  completely antisymmetric contravariant tensor of maximal degree with 
  component $\varepsilon^{1\dots D}=1$ in the canonical basis of 
  $\mathbb R^D$. Then the corresponding isomorphism $\ast:\wedge_{N} 
  \mathbb R^{D\ast} \rightarrow \wedge_{N}\mathbb R^D$ extends by 
  $C^\infty(\mathbb R^D)$-linearity as an isomorphism of $C^\infty 
  (\mathbb R^D)$-modules, again denoted by $\ast$, of 
  $\Omega_{N}(\mathbb R^D)$ into the space (of contravariant tensor 
  fields on $\mathbb R^D$) $\wedge_{N}\mathbb R^D\otimes C^\infty 
  (\mathbb R^D)$ with
  \[
  \ast \Omega^p_{N}(\mathbb R^D)= \wedge^{(N-1)D-p}_{N}\mathbb R^D 
  \otimes C^\infty (\mathbb R^D)
  \]
  for any $0\leq p\leq (N-1)D$. Let us define the first-order 
  differential operator $\delta$ of degree $-1$ on $\wedge_{N}\mathbb 
  R^D\otimes C^\infty (\mathbb R^D)$
  \[
  \delta : \wedge^{(N-1)p+r}_{N} \mathbb R^D\otimes C^\infty(\mathbb 
  R^D)\rightarrow \wedge^{(N-1)p+r-1}_{N}\mathbb R^D \otimes C^\infty 
  (\mathbb R^D)
  \]
  by setting
  
  \beq
  \delta T = \Yg^N_{(N-1)p+r-1}\circ \tilde \delta T
  \label{eq14}
  \end{equation}
   for $T\in \wedge^{(N-1)p+r}_{N}\mathbb R^D\otimes C^\infty(\mathbb 
    R^D)$ with $0\leq p<D$ and $1\leq r\leq N-1$,  $\tilde 
    \delta$ being defined by
  
  \[
  (\tilde\delta T)^
  {\mu^1_{1}\dots \mu^{p+1}_{1},\dots,\mu^1_{r-1}\dots 
    \mu^{p+1}_{r-1},\mu^1_{r}\dots 
    \mu^p_{r},\dots,\mu^1_{N-1}\dots \mu^p_{N-1}} =
    \partial_{\mu} T^{\mu^1_{1}\dots \mu^{p+1}_{1},\dots,\mu^1_{r}\dots 
    \mu^{p}_{r}\mu,\dots, 
    \mu^1_{N-1}\dots \mu^p_{N-1}}
    \]
    where we have used the canonical identification (\ref{Id}) and 
    the conventions explained below (\ref{Id}). It is worth noticing 
    here that in view (essentially) of Proposition 1, one has $\delta 
    T=\tilde \delta T$ for $r=N-1$, i.e. in this case (well-filled 
    case) the projection is not necessary in formula (\ref{eq14}).
    So defined $(\delta T)(x)$ is by construction in 
    $\wedge^{(N-1)p+r-1}_{N}\mathbb R^D$ and the operator 
    $\delta$ is in each degree proportional to the operator $\ast 
    d\ast^{-1}$, i.e. that one has
    \beq
    \delta = c_{n} \ast d\>\> \ast^{-1}:\wedge^n_{N}\mathbb R^D\otimes 
    C^\infty (\mathbb R^D) \rightarrow \wedge^{n-1}_{N}\mathbb 
    R^D\otimes C^\infty (\mathbb R^D)
    \label{eq15}
    \end{equation}
    for some $c_{n} \in \mathbb R$, $1\leq n\leq (N-1)D$ ($\delta=0$ in 
    degree zero).

\section{Digression on a related multicomplex} 

In this section, we 
introduce a  multicomplex  which will be related to our $N$-complex 
$\Omega_{N}(\mathbb R^D)$ in the next section. We also derive some 
useful cohomological results in this multicomplex, which will be the key 
for proving our generalization of the Poincar\'e lemma that is 
Theorem~\ref{theo1}.\\

Let $\fraca$ be the graded
tensor product of $N-1$ copies of the
exterior algebra $\wedge \mathbb R^{D\ast}$ of the dual space 
$\mathbb R^{D\ast}$ of $\mathbb R^D$ with $C^\infty( \mathbb R^D)$,
$$\fraca  = (\otimes^{N-1} \wedge \mathbb R^{D\ast}) \otimes C^\infty(\mathbb 
R^D)=\otimes^{N-1}_{C^\infty(\mathbb R^D)}\Omega(\mathbb R^D) . $$
An element of $\fraca$ is  as a sum of products of the $(N-1)D$ 
generators $d_i x^\mu$ ($i = 1, \dots, N-1$, $\mu = 1, \dots, D$) with 
 smooth functions on $\mathbb R^D$. Elements of $\fraca$ will be 
 refered to as {\sl multiforms}.
 The space $\fraca$ is a 
 graded-commutative algebra for the total degree, in particular one has 
$$ d_i x^\mu \, d_j 
x^\nu = - d_j x^\nu \, d_i x^\mu, \; \; \; \; x^\mu \, d_i x^\nu = d_i 
x^\nu \, x^\mu. $$ 
One defines $N-1$ antiderivations $d_i$ on $\fraca$ by 
setting
\begin{equation} 
 d_i f = \partial_\mu f \, d_i x^\mu \; \; (f \in 
C^\infty(\mathbb R^D))\, , \; \; \; \; \; \; d_i(d_j x^\mu) = 0. 
\end{equation}
These antiderivations anticommute,
\begin{equation}
d_i d_j + d_j d_i = 0
\end{equation}
in particular each $d_i$ is a differential. The graded algebra 
$\fraca$ has a natural multidegree  $(\fracd_1, \fracd_2, \dots, 
\fracd_{N-1})$ for which $\fracd_i(d_j x^\mu) = \delta_{ij}$.

It is useful to consider
the subspaces $\fraca^{(k)}$ of multiforms that vanish at the 
origin, together with all their successive derivatives up to order $k-1$ 
included ($k \geq 1$). If $\omega \in \fraca^{(k)}$, one says that 
$\omega$ is of {\sl order $k$}. The terminology
comes from the fact that a smooth function belongs to $\fraca^{(k)}$ if 
and only if
its limited Taylor expansion starts with terms of order $\geq k$.
If $l \geq k$, $\fraca^{(l)} \subset \fraca^{(k)}$. The subspaces 
$\fraca^{(k)}$ are
not stable under $d_i$ but one has $d_i \fraca^{(k)} \subset 
\fraca^{(k-1)}$ for $k \geq 1$ (with $\fraca^{(0)} \equiv \fraca$). The 
vector space $H^{(k)}(d_i, \fraca)$ is defined as
$$H^{(k)}(d_i, \fraca) \equiv \frac{Z^{(k)}(d_i, \fraca)} 
{d_i\fraca^{(k+1)}}$$ where $Z^{(k)}(d_i, \fraca)$ is the set of 
$d_i$-cocycles $\in \fraca^{(k)}$. Note that any multiform $\omega
\in \fraca$ can be written as $\omega = p^{(k)} + \beta$ where $p^{(k)}$ 
is a polynomial multiform of polynomial degree $k$ and $\beta \in \fraca^{(k+1)}$. 
This decomposition is unique which implies in 
particular that $\fraca^{(k)} \cap d_i\fraca
= d_i\fraca^{(k+1)}$.

It follows from the standard Poincar\'e lemma that 
\begin{equation}
H^{(1)}(d_i, \fraca) = 0.
\label{Poinc21}
\end{equation}
Indeed, the cohomology of $d_i$ in $\fraca$ is isomorphic to the space of 
constant multiforms not involving $d_i x^\mu$. The condition that the 
cocycles belong to $\fraca^{(1)}$, i.e., vanish at the origin, eliminates 
precisely the constants. One has also $H^{(m)}(d_i, \fraca) = 0 \; \forall 
m\geq 1$ since $\fraca^{(m)} \subset \fraca^{(1)}$ for $m \geq 1$ and 
$\fraca^{(m)} \cap d_i\fraca = d_i\fraca^{(m+1)}$. 

Let $K$ be an arbitrary subset of $\{1, 2, \dots, N-1\}$. We define 
$\fraca_K$ as
the quotient space
$$ \fraca_K = \frac{\fraca}{\sum_{j \in K} d_j \fraca} $$ 
(for $K = \emptyset$, $\fraca_K = \fraca$). The differential $d_i$ induces, for each 
$i$, a well-defined differential in $\fraca_K$ which we still denote by 
$d_i$. Of course, the induced $d_i$ is equal to zero if $i \in K$. 
\begin{lemma} For every proper subset $K$ of $\{1, 2, \dots, N-1\}$ and 
for every $i \notin K$, one has
$$ H^{(k+1)}(d_i , \fraca_K) = 0 \; \; \; 
(k = \# K)$$
\end{lemma}

\noindent {\bf Proof} The proof proceeds by induction on the
number $k$ of elements of $K$. The lemma clearly holds for $k=0$ ($K = 
\emptyset$) since then $\fraca_K = \fraca$ and the lemma reduces to 
Eq. (\ref{Poinc21}).
Let us now assume that the lemma holds for all subsets $K$ (not containing 
$i$) with $k \leq \ell$ elements. Let $K'$ be a subset not containing $i$ 
with $\ell+1$ elements. Let $j \in K'$ and $K'' = K' \backslash \{j\}$.
The induction hypothesis implies $H^{(\ell+1)}(d_i , \fraca_{K''}) = 
H^{(\ell+1)}(d_j , \fraca_{K''}) = 0$.
By standard
``descent equation" arguments (see below), this leads to $$ H^{p,q,(\ell+2)} (d_i \vert 
d_j , \fraca_{K''}) \simeq H^{p+1, q-1, (\ell+2)} (d_i \vert d_j , 
\fraca_{K''}).$$ In $H^{p,q,(\ell+2)} (d_i \vert d_j , \fraca_{K''})$, the 
first superscript $p$ stands for the $\fracd_i$-degree, the second supercript 
$q$ stands for the $\fracd_j$-degree while $(\ell+2)$ is the polynomial 
order.
Repeated application of this isomorphism yields
$$ H^{p,q,(\ell+2)} (d_i \vert d_j , \fraca_{K''}) \simeq 
H^{p+q,0,(\ell+2)} 
(d_i \vert d_j , \fraca_{K''}).$$ But $H^{p+q,0,(\ell+2)} (d_i \vert d_j , 
\fraca_{K''}) \equiv H^{p+q,0,(\ell+2)} (d_i , \fraca_{K''}) = 0$. Hence, the 
cohomological
spaces $H^{p,q,(\ell+2)} (d_i \vert d_j , \fraca_{K''})$ vanish for all $p$, 
$q$, which is precisely the statement $ H^{(\ell+2)} (d_i, \fraca_{K'}) = 
0$.$\square$ \\

The precise descent equation
argument involved in this proof runs as follows: let $\alpha^{p,q,(\ell+2)}$ be a $d_i$-cocycle 
modulo $d_j$ in $\fraca_{K''}$, i.e., a solution of $d_i 
\alpha^{p,q,(\ell+2)}
+ d_j \alpha^{p+1, q-1,(\ell+2)} \approx 0$ for some $\alpha^{p+1, 
q-1,(\ell+2)}$, where the notation $\approx$ means ``modulo terms in $\sum_{j 
\in K''} d_j \fraca$. Applying $d_i$ to this equation yields $d_j d_i 
\alpha^{p+1, q-1,(\ell+2)} \approx 0$ and hence, since $d_i \alpha^{p+1, 
q-1,(\ell+2)}$ is of order $\ell+1$ and $H^{(\ell+1)}(d_j , \fraca_{K''}) = 0$, 
$d_i \alpha^{p+1, q-1,(\ell+2)} + d_j \alpha^{p+2,q-2,(\ell+2)} \approx 0$ for 
some $\alpha^{p+2,q-2,(\ell+2)}$. Hence, $\alpha^{p+1, q-1,(\ell+2)}$ is also a 
$d_i$-cocycle modulo $d_j$ in $\fraca_{K''}$. Consider the map 
$\alpha^{p,q,(\ell+2)} \mapsto \alpha^{p+1, q-1,(\ell+2)}$ of $d_i$-cocycles 
modulo $d_j$. There is an arbitrariness in the choice of 
$\alpha^{p+1,q-1, (\ell+2)}$ given $\alpha^{p,q,(\ell+2)}$ so  this map 
is ambiguous, however $H^{(\ell+1)}(d_j , \fraca_{K''}) = 0$ implies 
that it induces a well-defined linear mapping
$H^{p,q,(\ell+2)} (d_i \vert d_j , \fraca_{K''}) \rightarrow H^{p+1, q-1, 
(\ell+2)} (d_i \vert d_j , \fraca_{K''})$ in cohomology. This map is 
injective and surjective since $H^{(\ell+1)}(d_i , \fraca_{K''}) = 0$ 
and thus one has the isomorphism 
$H^{p,q,(\ell+2)}(d_{i}\vert d_{j},\fraca_{K''})\simeq 
H^{p+1,q-1,(\ell+2)}(d_{j}\vert d_{i},\fraca_{K''})$  (see \cite{D-VHTV} 
for additional information).\\

A direct application of this lemma is the following 
\begin{proposition}
Let $J$ be any  non-empty
subset of $\{1,2, \dots, N-1\}$.
Then
$$\big(\prod_{j \in J} d_j \big) \alpha = 0 \hbox{ and } \alpha \in \fraca^{(\# J)}
\Rightarrow \alpha = \sum_{j \in J} d_j \beta_j $$
for some $\beta_j$'s. 
\end{proposition}

\noindent {\bf Proof} The property is clearly true for $\# J = 1$ (see Eq. 
(\ref{Poinc21})). Assume then that the property is true for all proper 
subsets with $k \leq \ell < N-1$ elements. Let $J$ be a proper subset with 
exactly $\ell$ elements and $i \notin J$. Let $\alpha$ be a multiform in 
$\fraca^{(\ell+1)}$ such that $d_i (\prod_{j \in J}d_{j}) \alpha = 0$. This is 
equivalent to $(\prod_{j \in J}d_{j}) d_i \alpha = 0$. Application of the 
recursive assumption to $d_i \alpha$, which belongs to $\fraca^{(\ell)}$, 
implies then $d_i \alpha = \sum_{j \in J}d_{j}\beta_j$, from which one 
derives, using the previous lemma, that $ \alpha = d_i \rho + \sum_{j \in 
J} \rho_j$ for some $\rho$, $\rho_j$. Therefore, the property passes on to 
all subsets with $\ell + 1$ elements, which establishes the 
theorem.$\square$\\

We are now in position to state and prove the main result of this 
section.

\begin{theorem}
Let $K$ be an arbitrary non-empty
subset of $\{1,2, \dots, N-1\}$. If the multiform $\omega$ is such that
\begin{equation}
\big( \prod_{i \in I} d_i \big) \omega = 0 \; \; \; \; \forall I \subset K 
\; \vert \, \# I = m \label{keykey1}
\end{equation}
(with $m \leq \# K$ a fixed integer), then 
\begin{eqnarray}
\omega =
\! \! \! \! \! \! \! \! \!
\sum_{\begin{array}{c}
J \subset K \\
\#J = \#K - m +1
\end{array}}
\! \! \! \! \! \! \! \! \!
\big( \prod_{j \in J} d_j \big) \alpha_J + \omega_0 \label{keykey2}
\end{eqnarray}
where $\omega_0$ is a polynomial multiform of degree $\leq m-1$.
\end{theorem}

\noindent {\bf Proof} The polynomial multiform $\omega_0 $ is clearly a solution of 
the problem, so we only need to check that if $\omega\in\fraca^{(m)}$ in 
addition to (\ref{keykey1}), then (\ref{keykey2}) is replaced by
\begin{eqnarray}
\omega =
\! \! \! \! \! \! \! \! \!
\sum_{\begin{array}{c}
J \subset K \\
\#J = \#K - m +1
\end{array}}
\! \! \! \! \! \! \! \! \!
\big( \prod_{j \in J} d_j \big) \alpha_J. \label{keykey3}
\end{eqnarray}
The $\alpha_J$'s can be assumed to be of order $\# K + 1$ since one 
differentiates them $\# K - m + 1$ times to get $\omega$. To prove 
(\ref{keykey3}), we proceed recursively, keeping $m$ fixed and increasing 
the size of $K$ step by step from $\# K = m$ to $\# K = N-1$. If $\# K = 
m$, there is nothing to be proven since $I = K$ and the theorem reduces to 
the previous theorem. So, let us assume that the theorem has been proven 
for $\# K = k \geq m$ and let us show that it extends to any set $U = K 
\cup \{\ell\}$, $\ell \notin K$ with $\# U = k+1$ elements. If (\ref{keykey1}) 
holds for any subset $I \subset U$ of $U$ (with $\# I = m$), it also holds 
for any subset
$I \subset K$ of $K \subset U$ (with $\# I = m$), so the recursive 
hypothesis implies
\begin{eqnarray}
\omega =
\! \! \! \! \! \! \! \! \!
\sum_{\begin{array}{c}
J \subset K \\
\#J = k - m +1
\end{array}}
\! \! \! \! \! \! \! \! \!
\big( \prod_{j \in J} d_j \big) \alpha_J. 
\label{recursive21}
\end{eqnarray}
Let now $A$ be an arbitrary subset of $U$ with $ \# A = m$, which contains 
the added element $\ell$. Among all the subsets $J$ occurring in the sum 
(\ref{recursive21}), there is only one, namely $J' = U \backslash A$ such 
that $J' \cap A = \emptyset$. The condition (\ref{keykey1}) with $I = A$ 
implies, when applied to the expression (\ref{recursive21}) of $\omega$, 
$$ \big( \prod_{j \in U} d_j \big) \alpha_{J'} = 0 $$ (if $J \not= J'$, 
the product $(\prod_{i \in A} d_i) (\prod_{j \in J}d_j
)$ identically vanishes because at least one differential $d_f$ is 
repeated). But since $\alpha_{J'}$ is of order $k+1 = \# U$, the previous 
proposition implies that $$ \alpha_{J'} = \sum_{j \in U} d_j \beta_{j} 
^{J'}.$$ When injected into (\ref{recursive21}), this yields in turn 
\begin{eqnarray}
\omega =
\! \! \! \! \! \! \! \! \!
\sum_{\begin{array}{c}
L \subset U \\
\#L = k - m +2
\end{array}}
\! \! \! \! \! \! \! \! \!
\big( \prod_{j \in L} d_j \big) \alpha'_L. 
\end{eqnarray}
for some $\alpha'_L$,
and shows that the required property is also valid for sets with cardinal 
equal to $k+1$,
completing the proof of the theorem.$\square$
 
\section{The generalization of the Poincar\'e lemma}

With the result of last section, Theorem 2, we can now proceed to the 
proof of Theorem 1 that is to the proof of the generalization of the 
Poincar\'e lemma announced in the introduction.\\

Let us first show that $\Omega_{N}(\mathbb R^D)$ identifies 
canonically as graded $C^\infty(\mathbb R^D)$-module with the image of 
a $C^\infty(\mathbb R^D)$-linear homogeneous projection $\pi$ of 
$\fraca$ into itself: $\Omega_{N}(\mathbb R^D)=\pi(\fraca)\subset 
\fraca$. Indeed by using the canonical identification (\ref{Id}) of 
Section 3, one has the identification
\begin{equation}\label{Sub}
 \wedge^{(N-1)n+i}_{N}E \subset \underbrace{\wedge^{n+1} E\otimes 
 \dots \otimes \wedge^{n+1}E}_{i}\otimes 
 \underbrace{\wedge^nE\otimes \dots \otimes \wedge^n E}_{N-1-i}
 \end{equation}
 of the Schur module $E^{Y^N_{(N-1)n+i}}=\wedge^{(n+1)n+i}_{N}E$ as 
 vector subspace of the right-hand side. However by decomposing the 
 right-hand side of (\ref{Sub}) into irreducible subspaces for the 
 action of $GL(E)$ one sees that there is only one irreducible factor 
 isomorphic to $E^{Y^N_{(N-1)n+i}}$ which is therefore the image of a 
 $GL(E)$-invariant projection. It follows that $\wedge_{N}E\subset 
 \otimes^{N-1} \wedge E$ is the image of a $GL(E)$-invariant 
 projection $P$ of $\otimes^{N-1}\wedge E$ into itself which is 
 homogeneous for the total degree. The result for $\Omega_{N}(\mathbb 
 R^D)$ follows by chosing $E$ to be the dual space $\mathbb R^{D\ast}$ of  
 $\mathbb R^D$ and by setting $\pi=P\otimes I_{C^\infty(\mathbb 
 R^D)}$ in view of $\Omega_{N}(\mathbb R^D)= \wedge_{N}\mathbb 
 R^{D\ast}\otimes C^\infty (\mathbb R^D)$ and $\fraca 
 =(\otimes^{N-1}\wedge \mathbb R^{D\ast})\otimes C^\infty(\mathbb 
 R^D)$. The projection $\pi$ is in fact by construction a projection 
 of $\oplus_{p\in \mathbb N}\fraca^{[p]}$ into itself where 
 $\fraca^{[p]}=\fraca^{n+1,\dots,n+1,n,\dots,n}$, $p=(N-1)n+i$ with 
 obvious notations.\\
 
 We now relate the $N$-differential $d$ of $\Omega_{N}(\mathbb R^D)$ 
 to the differentials $d_{i}$ of $\fraca$. Let $\omega$ be an element 
 of $\Omega^{p}_{N}(\mathbb R^D)$ with $p=(N-1)n+i$,  $0\leq i< N-1$. 
 One has

\begin{equation}
d\omega = c_\omega \pi(d_{i+1} \omega)
\end{equation}
where $c_\omega$ is a non-vanishing number that depends
on the degrees of $\omega$.
In general, the projection is non trivial, in
the sense that $d_{i+1} \omega$ has components not only along
the irreducible Schur module $E^{Y^N_{p+1}}$ ($E=\mathbb R^{D\ast}$), but also along other
Schur modules not occurring in $\Omega_N(\mathbb R^D)$.
For instance, with $N=3$, the covariant vector with components $v_\alpha$ defines
the element $v = v_\alpha d_1 x^\alpha$ of $\fraca$.
One has $d_2 v =  - \partial_\beta v_\alpha d_1x^\alpha d_2
x^\beta$.  This expression contains both a symmetric ($dv$) {\em and} an
antisymmetric part, so 
$d_2 v = d v -  v_{[\alpha, \beta]} d_1 x^\alpha d_2
x^\beta$.  The projection removes $v_{[\alpha, \beta]} d_1 x^\alpha d_2
x^\beta$, which does not vanish in general.  
Because the projection is non-trivial, the conditions
$d \omega = 0$ and $d_{i+1} \omega = 0$ are inequivalent for generic $i$.
However, if $\omega$ is a well-filled tensor that is if $i = 0$, then 
\begin{equation}
d \omega = d_1 \omega \; \; \; (i=0)
\end{equation}
Indeed, $d_1 \omega$ has automatically the correct Young symmetry.
Thus the conditions
$d_1 \omega = 0$ and $d \omega = 0$ are equivalent.
Furthermore, because of the symmetry between the columns, if 
$d_1 \omega = 0$, then, one has also $d_2 \omega =
d_3 \omega = \cdots = 0$.  For instance, again
for $N=3$, the derivative of the symmetric tensor
$g = g_{\alpha \beta} d_1 x^\alpha d_2 x^\beta$
($g_{\alpha \beta} = g_{\beta \alpha}$) is given by
$dg = d_1 g = \frac{1}{2} (g_{\alpha \beta , \rho}
- g_{\rho \beta , \alpha}) d_1 x^\rho d_1 x^\alpha d_2 x^\beta$.
The completely symmetric component $g_{(\alpha \beta, \rho)}$
is absent because $d_1 x^\rho d_1 x^\alpha = -
d_1 x^\alpha d_1 x^\rho$.  Also, it is clear that if
$d_1 g = 0$, then, $d_2 g = \frac{1}{2} (g_{\alpha \beta , \rho}
- g_{\alpha \rho , \beta}) d_1 x^\alpha d_2 x^\beta d_2 x^\rho
= 0$.
This generalizes to the following lemma:
\begin{lemma}
Let $\omega \in \Omega_N^{(N-1)n}(\mathbb R^D)$
(well-filled, or rectangular, tensor).  Then
\begin{equation}
d^k \omega = 0  \; \; \; \; \;
\Leftrightarrow \; \; \; \; \; \;
(\prod_{j \in J, \, \# J = k} \, d_j)
\omega = 0.
\end{equation}
\end{lemma}

{\bf Proof} One has $d^k \omega = (-1)^m d_1 d_2 \cdots d_k \omega$.
Indeed, it is clear that the multiform
 $ d_1 d_2 \cdots d_k \omega \in \fraca^{n+1,
n+1, \cdots, n+1, n , \cdots, n}$ belongs  to
$\Omega_N(\mathbb R^D)$ because it cannot have components
along the invariant subspaces corresponding to Young diagrams
with first column having $i >r + 1$ boxes, since
one cannot put two derivatives $\partial_\mu$, $\partial_\nu$
in the same column.  Hence, $d^k \omega = 0$ is equivalent to
$d_1 d_2 \cdots d_k \omega = 0$.  One completes the proof by
observing that for well-filled tensors, the
condition $d_1 d_2 \cdots d_k \omega = 0$ is equivalent to
the conditions $d_{i_1} d_{i_2} \cdots d_{i_k} \omega
= 0 \; \; \forall (i_1, \cdots, i_k)$ because of symmetry in
the columns.
$\square$\\

\begin{lemma}
Let $\omega \in \Omega_N^{(N-1)n}(\mathbb R^D)$ with $n\geq 1$.  Then
\begin{equation}
\omega = \sum_{J,\, \# J = N-k} \big( \prod_{j \in J}
d_j \big) \, \alpha_J \; \; \Rightarrow
\omega = d^{N-k} \alpha
\label{implic21}
\end{equation}
for some $\alpha \in \Omega_N^{(N-1)n-N+k}(\mathbb R^D)$, $k\in 
\{1,\dots,N-1\}$.
\end{lemma}
{\bf Proof} First, we note that the $\alpha_J$ occurring
in (\ref{implic21}) can be chosen to
have $\fracd_{i}$-degrees  equal to $n-1$ or $n$
according to whether $d_i$ acts or does not act on
$\alpha_J$, since $\omega$
has multidegree $(n,n,\cdots,n)$.
Second, one can project the right-hand side of (\ref{implic21})
on $\Omega_N^{(N-1)n}(\mathbb R^D)$ without changing the
left-hand side, since $\omega \in \Omega_N^{(N-1)n}(\mathbb R^D)$.
It is easy to see that
$\pi [( \prod_{j \in J}
d_j ) \, \alpha_J] \sim d^{N-k} \alpha'_J$, with $\alpha'_J=
\pi( \tilde{\alpha}_J)$, where  $\tilde{\alpha}_J$ is the element
in $\fraca^{n,\cdots, n, n-1, n-1, \cdots, n-1}$
obtained by reordering the ``columns" of $\alpha_J$
so that they have non-increasing
length. 
In fact, when differentiated,
the other irreducible components of
$ \tilde{\alpha}_J$ do not contribute to
$\omega$ because their first
column is too long to start with or because two partial
derivatives find themselves in the same column, yielding zero.
Injecting the above  expression for $\pi [( \prod_{j \in J}
d_j ) \, \alpha_J]$ into (\ref{implic21}) yields the desired 
result.~$\square$\\

\begin{lemma}  Let $\omega \in \Omega_N^{(N-1)n}(\mathbb R^D)$ 
with $n \geq 1$ be
a polynomial multiform of degree $k-1$.  Then, 
\begin{equation}
\omega = d^{N-k} \alpha
\end{equation}
for some polynomial multiform $\alpha \in \Omega_N^{(N-1)n-N+k}
(\mathbb R^D)$
of degree $N-1$, with $k\in \{1,\dots,N-1\}$.
\end{lemma}

{\bf Proof} The proof amounts to play with Young diagrams.  
The coefficients
of $\omega$ transform in the tensor product 
of the representation associated with $Y^N_{(N-1)n}$
(symmetry of $\omega$) and the completely symmetric
representation with $k-1$ boxes (symmetric polynomials
in the $x^\mu$'s of degree $k-1$).   Let $T$ be
this representation and $V_T$ be the carrier vector space.
Similarly, the multiform
$\alpha$ transforms (if it exists) in the
tensor product of the representation associated
with  $Y^N_{(N-1)n-N+k}$ (symmetry of $\alpha$) and the completely symmetric
representation with $N-1$ boxes (symmetric polynomials
in the $x^\mu$'s of degree $N-1$).  Let $S$ be this
representation and $W_S$ be the carrier vector space.
Now, the linear
operator $d^{N-k}: W_S \rightarrow V_T$ is an intertwiner 
for the representations
$S$ and $T$.  To analyse how it acts, it is convenient to
decompose both $S$ and $T$ into irreducible representations.

The crucial fact is that all irreducible
representations occurring in
$T$ also occur in $S$.  That is, if 
$$
T = \oplus_i T_i, \; \;  \; \; V_T = \oplus_i V_i
$$
(where each irreducible representation $T_i$ has 
multiplicity one), then
$$
S = (\oplus_i T_i) \oplus (\oplus_\alpha T_\alpha), \; \; \; \;
W_S = (\oplus_i W_i) \oplus (\oplus_\alpha W_\alpha)
$$
where $T_\alpha$ are some other representations, irrelevant
for our purposes. Because $T_i$ is irreducible,
the operator $d^{N-k}$ maps the invariant
subspace $W_i$ on the
invariant subspace $V_i$, and furthermore, $d^{N-k}\vert_{W_i}$
is either zero or bijective.  It is easy to verify by taking simple
examples that $d^{N-k}\vert_{W_i}$ is not zero.  Hence,
$d^{N-k}\vert_{W_i}$ is injective, which implies
that $d^{N-k}: W_S \rightarrow V_T$ is surjective, so that
$\omega$ can indeed be written as $d^{N-k} \alpha$ for some $\alpha$. 
~$\square$\\

{\bf Proof of Theorem \ref{theo1}} The theorem \ref{theo1}
is a direct consequence of the above two lemmas.  (i) Let
$\omega \in \Omega_N^{(N-1)n}(\mathbb R^D)$ (with $n \geq 1$)
be annihilated
by $d^k$, $d^k \omega = 0$. We write $\omega = \omega' +
\omega_0$, where $\omega'$ is of order $k$  and
where $\omega_0$  is a polynomial multiform of polynomial degree
$<k$.  Both $\omega'$ and $\omega_0$ have the symmetry of
$\omega$.  Also, since $\omega_0$ is trivially annihilated by $d^k$, one
has separately $d^k \omega' = 0$ and $d^k \omega_0
=0$.
We consider first $\omega'$.  The first lemma implies
$(\prod_{j \in J, \, \# J = k} \, d_j)
\omega' = 0$, from which it follows, using the theorem of the previous
section, that $$\omega'  = \sum_{J,\, \# J = N-k} \big( \prod_{j \in J}
d_j \big) \, \alpha_J $$ (see (\ref{keykey3})).
By the
second lemma above, this term can be written as $d^{N-k} \alpha$.
As we have also seen, the same property holds for $\omega_0$.
This proves the theorem for $n \geq 1$.
(ii) The case of $H^0_{(k)} (\Omega_N(\mathbb R^D))$
is even easier to discuss: for a function, the condition $d^k f = 0$
is equivalent to $\partial_{\mu_1 \cdots \mu_k} f = 0$ and
thus, $f$ must be of degree strictly less than $k$.
Moreover, it can never be the $d^{N-k}$ of something, since
there is nothing in negative degree.  $\square$\\

It is worth noticing here that, as explained in the 
introduction, Theorem~\ref{theo1} has a dual counterpart for the 
$\delta$-operator introduced at the end of Section 3 which allows to 
integrate lots of generalized currents conservation equations.\\
In the last section of this paper we shall sketch another approach 
for proving Theorem~\ref{theo1} which is based on the appropriate
generalization of homotopy for $N$-complexes.

\section{Structure of $H^m_{(k)}(\Omega_{N}(\mathbb R^D))$ for 
generic $m$}

In the previous section we have shown that 
$H^m_{(k)}(\Omega_{N}(\mathbb R^D))$ vanishes whenever $m=(N-1)n$ with 
$n\geq 1$. In the case $N=2$ this is the usual Poincar\'e lemma which 
means that the cohomology vanishes in positive degrees. For $N\geq 3$ 
there are degrees $m$ which do not belong to the set $\{ 
(N-1)(n+1)\vert n\in \mathbb N\}$ and it turns out that for such a 
(generic) degree $m$, the spaces $H^m_{(k)}(\Omega_{N}(\mathbb R^D))$ 
are non trivial $(k\in \{1,\dots, N-1\})$. More precisely for $m\in 
\{0,\dots, N-2\}$ these spaces are finite-dimensional of strictly 
positive dimensions whereas for $m\geq N$ with $m\not= (N-1)n$ these 
spaces are infinite-dimensional. In the following we shall 
explicitly display the case $N=3$ and indicate how to proceed for 
the general case $N\geq 3$.\\

In order to simplify the notations let us denote the spaces 
$H^m_{(k)}(\Omega_{N}(\mathbb R^D))$ by $H^m_{(k)}$ and the graded 
spaces $H_{(k)}(\Omega_{N}(\mathbb R^D))$ by $H_{(k)} 
(=\oplusinf_{m}H^m_{(k)})$. \\

For $N=3$, one has only $H_{(1)}$ and $H_{(2)}$ and Theorem 1 states 
that $H^{2n}_{(1)}=H^{2n}_{(2)}=0$ for $n\geq 1$ and that 
$H^0_{(1)}\simeq \mathbb R$ is the space of constant functions on 
$\mathbb R^D$ whereas $H^0_{(2)}$ is the space of polynomial functions 
of degree less or equal to one on $\mathbb R^D$, i.e. 
$H^0_{(1)}\simeq \mathbb R \oplus \mathbb R^D$. On the other hand, the 
elements of $H^1_{(1)}$ identify with the covariant vector fields (or 
one-forms) $x\mapsto X(x)$ on $\mathbb R^D$ satisfying
\begin{equation}\label{Killing}
 \partial_{\mu} X_{\nu} +\partial_{\nu} X_{\mu} =0
 \end{equation}
 which is the equation characterizing the Killing vector fields (i.e. 
 infinitesimal isometries) of the standard euclidean metric 
 $\sum^D_{\mu=0}(dx^\mu)^2$ of $\mathbb R^D$. The general solution of 
 (\ref{Killing}) is $X_{\mu}=v_{\mu}+a_{\mu\nu}x^\nu$ with $v\in 
 \mathbb R^D$ (infinitesimal translations) and $a\in \wedge^2\mathbb 
 R^D$ i.e. $a_{\mu\nu}=-a_{\nu\mu}$=C$^{\mbox{te}}$ 
 (infinitesimal rotations). Thus one has $H^1_{(1)}\simeq \mathbb 
 R^D\oplus \wedge^2 \mathbb R^D$. Notice that with this terminology we 
 have implicitly identified covariant vector fields with 
 contravariant ones by using the standard metric of $\mathbb R^D$. 
 Notice also that as far as $H^0_{(1)}$, $H^0_{(2)}$ and $H^1_{(1)}$ 
 are concerned nothing change if $N\geq 3$. For $N=3$, $H^1_{(2)}$ 
 identifies with the space of covariant vector fields $x\mapsto  
 X(x)$ on $\mathbb R^D$ satisfying
 \begin{equation}
  \partial_{\lambda}(\partial_{\mu}X_{\nu}-\partial_{\nu}X_{\mu})=0
  \label{eq.6.2}
 \end{equation}
 modulo the ones of the form $X_{\mu}=\partial_{\mu}\varphi$ for 
 some $\varphi\in C^\infty(\mathbb R^D)$. The general solution of 
 (\ref{eq.6.2}) is 
 $X_{\mu}=a_{\mu\nu}x^\nu+\partial_{\mu}\varphi$ with 
 $a \in \wedge^2\mathbb R^D$ and $\varphi\in C^\infty (\mathbb 
 R^D)$ so that one has $H^1_{(2)}\simeq \wedge^2\mathbb R^D$. Let us 
 now show that $H^3_{(1)}$ is infinite-dimensional for $N=3$. For 
 this, consider an arbitrary 2-form $\omega$ i.e. an arbitrary 
 covariant antisymmetric tensor field of degree 2 on $\mathbb R^D$ and 
 consider the element $t={\mathbf Y}^3_{3}\circ 
 \stackrel{(0)}{\nabla}\omega$ of $\Omega^3_{3}(\mathbb R^D)$. Up to 
 an irrelevant normalization constant, the components of $t$ are given 
 by
 \begin{equation}
t_{\mu\lambda\nu}=2\partial_{\lambda}\omega_{\mu\nu}+\partial_{\mu}
  \omega_{\lambda\nu}-\partial_{\nu}\omega_{\lambda\mu} 
  \label{eq.6.3}
  \end{equation}
  and one verifies that one has $dt=0$ in $\Omega_{3}(\mathbb R^D)$. 
  On the other hand one has $t=dh$ in $\Omega_{3}(\mathbb R^D)$ that 
  is
  \begin{equation}
 2\partial_{\lambda}\omega_{\mu\nu}+\partial_{\mu}\omega_{\lambda\nu}- 
 \partial_{\nu}\omega_{\lambda\mu}=\partial_{\nu}
 h_{\mu\lambda}-\partial_{\mu}h_{\nu\lambda}
 \label{eq.6.4}
 \end{equation}
 for some symmetric covariant tensor field $h\in \Omega^2_{3}(\mathbb R^D)$ 
 if and only if $\omega$ is of the form
 \begin{equation}
  \omega_{\mu\nu}=a_{\rho\mu\nu}x^\rho + 
  \partial_{\mu}X_{\nu}-\partial_{\nu}X_{\mu}
 \label{eq.6.5}
 \end{equation}
 for $a\in \wedge^3\mathbb R^D$ and some covariant vector field $X\in \Omega^1_{3}(\mathbb R^D)$ 
 and then $t$ is proportional to $d^2(X)$ in $\Omega_{3}(\mathbb R^D)$ i.e. $t$ is trivial in $H^3_{(1)}$. This 
 argument shows firstly that $H^3_{(1)}$ contains the quotient 
 of the space of 2-forms by the ones of the form given by 
 (\ref{eq.6.5}) which is infinite-dimensional and secondly that the 
 same space identifies with a subspace of 
 $H^3_{(2)}$ which is therefore also infinite-dimensional. In fact as 
 will be shown below  one has an isomorphism
 $H^3_{(1)}\simeq H^3_{(2)}$ which is induced by the inclusion $\ker 
 (d)\subset \ker (d^2)$. By replacing the 2-form $\omega$ by an 
 irreducible covariant tensor field $\omega_{n}$ of degree $2n+2$ on 
 $\mathbb R^D$ with Young symmetry type given by the Young diagram 
 with $n$ lines of length two and two lines of length one, it can be 
 shown similarily that $H^{2(n+1)+1}_{(1)}$ and $H^{2(n+1)+1}_{(2)}$ 
 are infinite-dimensional spaces (we shall see that they are in fact isomorphic).\\
 
 The last argument for $N=3$ admits the following generalization for 
 $N\geq 3$. Let $Y^N_{m}$ be a Young diagram of the sequence 
 $(Y^N)$ and let $Y'_{m-p}$ be a Young diagram obtained by deleting 
 $p$ boxes of $Y^N_{m}$ with $0<p<N-1$ such that it does not belong to $(Y^N)$ 
 (i.e. $Y'_{m-p}\not= Y^N_{m-p})$ and such that by applying $p$ 
 derivatives (i.e. $\stackrel{(0)}{\nabla^p}$) to a generic tensor 
 field with Young symmetry $Y'_{m-p}$ one obtains a tensor which has 
 a nontrivial component $t$ with Young symmetry $Y^N_{m}$. Then 
 generically the latter $t\in \Omega^m_{N}(\mathbb R^D)$ is a 
 nontrivial generalized cocycle and one obtains by this procedure an 
 infinite dimensional subspace of the corresponding generalized 
 cohomology, i.e. of $H^m_{(k)}$ for the appropriate $k$. Notice that 
 this is only possible for $m\geq N$ with $m\not= (N-1)n$. We 
 conjecture that the whole nontrivial part of the generalized 
 cohomology of $\Omega_{N}(\mathbb R^D)$ in degree $m\geq N$ is 
 obtained by the above construction ($N\geq 3$).\\
 
 In order to complete the discussion for $N\geq 3$ in degree $m\leq 
 N-2$ as well as to show the isomorphisms $H^{2n+1}_{(1)}\simeq 
 H^{2n+1}_{(2)}$ for $N=3$, $n\geq 1$ and their generalizations for 
 $N\geq 3$, we now recall a basic lemma of the general theory of 
 $N$-complexes \cite{D-V2}, \cite{D-VK}. This lemma was formulated 
 in \cite{D-V2} in the more general framework of $N$-differential 
 modules (Lemma 1 of \cite{D-V2}) that is of ${\mathbf k}$-modules 
 equipped with an endomorphism $d$ such that $d^N=0$ where ${\mathbf k}$
 is a unital commutative ring. In this paper we only discuss $N$-complexes of 
 (real) vector spaces. Let $E$ be a $N$-complex of 
 cochain \cite{D-V2} like $\Omega_{N}(\mathbb R^D)$, that is here 
 $E=\oplusinf_{m\in \mathbb N}E^m$ is a graded vector space equipped 
 with an endomorphism $d$ of degree one such that $d^N=0$ ($N\geq 
 2$). The inclusions $\ker(d^k)\subset \ker(d^{k+1})$ and 
 $\im(d^{N-k})\subset \im (d^{N-k-1})$ induce  linear 
 mappings $[i]:H_{(k)}\rightarrow H_{(k+1)}$ in generalized cohomology 
 for $k$ such that $1\leq k\leq N-2$. Similarily the linear mappings 
 $d:\ker(d^{k+1})\rightarrow \ker(d^k)$ and 
 $d:\im(d^{N-k-1})\rightarrow \im(d^{N-k})$ obtained by restriction of 
 the $N$-differential $d$ induce linear mappings 
 $[d]:H_{(k+1)}\rightarrow H_{(k)}$. One has the following lemma (for 
 a proof we refer to \cite{D-VK} or \cite{D-V2}).
 
 \begin{lemma}\label{lemfond}
  
Let the integers $k$ and $\ell$  be such that $1\leq k$, $1\leq \ell$, 
$k+\ell\leq N-1$. Then the 
 hexagon of linear mappings
 \[ \begin{diagram} \node{}
\node{H_{(\ell+k)}(E)} \arrow{e,t}{[d]^k} \node{H_{(\ell)}(E)}
\arrow{se,t}{[i]^{N-(\ell+k)}} \node{} \\ \node{H_{(k)}(E)}
\arrow{ne,t}{[i]^\ell}
\node{} \node{} \node{H_{(N-k)}(E)} \arrow{sw,b}{[d]^\ell} \\ \node{}
\node[1]{H_{(N-\ell)}(E)} \arrow{nw,b}{[d]^{N-(\ell+k)}}
\node{H_{(N-(\ell+k))}(E)}
\arrow{w,b}{[i]^k} \node{} \end{diagram} \]
is exact.
\end{lemma}
Since $[i]$ is of degree zero while $[d]$ is of degree one, 
these hexagons give long exact sequences.\\

Let us apply the above result to the $N$-complex $\Omega_{N}(\mathbb 
R^D)$. For $N=3$, there is only one hexagon as above $(k=\ell=1)$ and, 
by using $H^{2n}_{(k)}=0$ for $n\geq 1$, $k=1,2$ it reduces to the 
exact sequences
\begin{equation}
 0\stackrel{[d]}{\rightarrow} H ^0_{(1)}\stackrel{[i]}{\rightarrow} 
 H^0_{(2)}\stackrel{[d]}{\rightarrow} 
 H^1_{(1)}\stackrel{[i]}{\rightarrow} H^1_{(2)}\stackrel{d}{\rightarrow} 
 0
 \label{eq.6.6}
 \end{equation}
 and 
 \begin{equation}
  0\stackrel{d}{\rightarrow} H^{2n+1}_{(1)}\stackrel{[i]}{\rightarrow} 
  H^{2n+1}_{(2)}\stackrel{d}{\rightarrow} 0
  \label{eq.6.7}
  \end{equation}
  for $n\geq 1$. The sequences (\ref{eq.6.7}) give the announced 
  isomorphisms $H^{2n+1}_{(1)}\simeq H^{2n+1}_{(2)}$ while the 
  4-terms sequence (\ref{eq.6.6}) allows to compute the finite 
  dimension of $H^1_{(2)}$ knowing the one of $H^0_{(1)}$, 
  $H^0_{(2)}$ and $H^1_{(1)}$. For $N\geq 3$ one has several hexagons 
  and by using $H^{(N-1)n}_{(k)}=0$ for $n\geq 0$, the sequence 
  (\ref{eq.6.6}) generalizes as the following $\frac{(N-2)(N-1)}{2}$ 
  four-terms exact sequences
  \begin{equation}
   0\stackrel{[d]^k}{\longrightarrow} 
   H^{k-1}_{(\ell)}\stackrel{[i]^{N-k-\ell}}{\longrightarrow} 
   H^{k-1}_{(N-k)}\stackrel{[d]^\ell}{\longrightarrow} 
   H^{k+\ell-1}_{(N-k-\ell)}\stackrel{[i]^k}{\longrightarrow }
   H^{k+\ell-1}_{(N-\ell)}\stackrel{[d]^{N-k-\ell}}{\longrightarrow} 0
   \label{eq.6.8}
   \end{equation}
   for $1\leq k,\ell$ and $k+\ell\leq N-1$. There are also two-terms 
   exact sequences generalizing (\ref{eq.6.7}) giving similar 
   isomorphisms but, for $N>3$, there are other longer exact 
   sequences (which are of finite lengths in view of 
   $H^{(N-1)n}_{(k)}=0$ for $n\geq 1$). Suppose that the spaces 
   $H^m_{(k)}$ are finite-dimensional for $k+m\leq N-1$ and that we 
   know their dimensions. Then the exact sequences (\ref{eq.6.8}) 
   imply that all the $H^m_{(k)}$ for $m\leq N-2$ are 
   finite-dimensional and allows to compute their dimensions in terms 
   of the dimensions of the $H^m_{(k)}$ for $k+m\leq N-1$. To complete 
   the discussion it thus remains to show the finite-dimensionality of 
   the $H^m_{(k)}$ for $k+m\leq N-1$. For $k+m\leq N-1$, the space 
   $H^m_{(k)}$ identifies with the space of (covariant) symmetric 
   tensor fields $S$ of degree $m$ on $\mathbb R^D$ such that
   \begin{equation}
    \sum_{\pi \in \cals_{k+m}}\partial_{\mu_{\pi(1)}}\dots 
    \partial_{\mu_{\pi(k)}} S_{\mu_{\pi(k+1)}}\dots \mu_{\pi(k+m)}=0
    \label{eq.6.9}
    \end{equation}
    for $\mu_{i}\in \{1,\dots,D\}$ where $\cals_{k+m}$ is the group of 
    permutation of $\{1,\dots,k+m\}$. In particular, for $k=1$ the 
    equation (\ref{eq.6.9}) means that $S$ is a Killing tensor field 
    of degree $m$ for the canonical metric of $\mathbb R^D$ and it is 
    well known and easy to show that the components of such a Killing 
    tensor field of degree $m$ are polynomial functions on $\mathbb 
    R^D$ of degree less or equal to $m$. In fact the Killing tensor 
    fields on $\mathbb R^D$ form an algebra for the symmetric product 
    over each point of $\mathbb R^D$ which is generated by the 
    Killing vector fields (which are polynomial of degree $\leq 1$). 
    Thus $H^m_{(1)}$ is finite-dimensional for $1+m\leq N-1$. By 
    using this together with Theorem 1, one shows by induction on $k$ 
    that $H^m_{(k)}$ is finite-dimensional for $k+m\leq N-1$, more 
    precisely, that the solutions of (\ref{eq.6.9}) are polynomial 
    functions on $\mathbb R^D$ of degree less than $k+m$.\\
    
    The results of this section concerning the generic degrees show 
    that our generalization of the Poincar\'e lemma, i.e. Theorem~\ref{theo1}, 
    is far from being a straightforward result and that it is optimal.

\section{Algebras}

Let $E\simeq \mathbb R^D$ be a $D$-dimensional vector space, $(Y)$ be 
a sequence $(Y)=(Y_{p})_{p\in \mathbb N}$ of Young diagrams such 
that $\vert Y_{p}\vert =p$ ($\forall p \in \mathbb N$) and let us use 
the notations and conventions of Section 3. As we have seen, the 
graded space $\wedge_{(Y)}E$ is a generalization of the exterior 
algebra of $E$ in the sense that as graded vector space it reduces to 
the latter when $(Y)$ coincides with the sequence 
$(Y^2)=(Y^2_{p})_{p\in\mathbb N}$ of the one-column Young diagrams. 
It is also a generalization of the symmetric algebra of $E$ since it 
reduces to it when $(Y)$ coincides with the sequences $(\tilde 
Y^2)=(\tilde Y^2_{p})_{p\in \mathbb N}$ of the one-line Young 
diagrams (which are the conjugates of the $Y^2_{p}$).
In fact, for general $(Y)$ the graded vector space
 $\wedge_{(Y)}E$ is also a graded algebra if one defines the product by 
 setting
 \begin{equation}
  TT'=\Yg_{p+p'} (T\otimes T')
  \label{eq8}
  \end{equation}
  for $T\in E^{Y_{p}}$ and $T'\in E^{Y_{p'}}$ where $\Yg_{n}$ is the 
  Young symmetrizer defined in Section 2. However, although it 
  generalizes the exterior product,  this product is 
  generically a nonassociative one. Thus $\wedge_{(Y)}E$ is a 
  generalization of the exterior algebra $\wedge E$ in which each 
  homogeneous subspace is irreducible for the action of $GL(E)\simeq 
  GL_{D}$ but in which one loses the associativity of the product. 
  There is another closely related generalization of the exterior 
  algebra connected with the sequence $(Y)$ in which what is retained 
  is the associativity of the graded product but in which one 
  generically loses the $GL(E)$-irreducibility of the homogeneous 
  components. This generalization, denoted by $\wedge_{[(Y)]}E$, is 
  such that $\wedge_{(Y)}E$ is a graded (right) 
  $\wedge_{[(Y)]}E$-module. We now describe its construction.

  Let $\Tg(E)$ be the tensor algebra 
  of $E$, we use the product defined by (\ref{eq8}) to equip 
  $\wedge_{(Y)}E$ with a right $\Tg(E)$-module structure by setting
  \begin{equation}
   T \lambda_{(Y)}(X_{1}\otimes \dots \otimes X_{n})=(\cdots(T 
   X_{1})\cdots )X_{n}
   \label{eq9}
   \end{equation}
   for any $X_{i}\in E$ and $T\in \wedge_{(Y)}E$. By definition the 
   kernel $\ker(\lambda_{(Y)})$ of $\lambda_{(Y)}$ is a two-sided ideal 
   of $\Tg(E)$ so that the right action of $\Tg(E)$ on $\wedge_{(Y)}E$ is 
   in fact an action of the quotient algebra  
   $\wedge_{[(Y)]}E=\Tg(E)/\ker(\lambda_{(Y)})$. So $\wedge_{(Y)}E$ is 
   a graded right $\wedge_{[(Y)]}E$-module.
   
   \begin{lemma}
    Let $N$ be an integer with $N\geq 2$ and assume that $(Y)$ is such 
    that the number of columns of the Young diagram $Y_{p}$ is 
    strictly smaller than $N$ for any $p\in \mathbb N$. Then 
    $\ker(\lambda_{(Y)})$ contains the two-sided ideal of $\Tg(E)$ which 
    consists of the tensors which are symmetric with respect to at 
    least $N$ of their entries; in particular 
    $(\lambda_{(Y)}(X))^N=0$, $\forall X\in E$.
    \end{lemma}
    Stated differently, under the assumption of the lemma for $(Y)$, 
    a monomial $X_{1}\dots X_{n}\in \wedge_{[(Y)]}E$ with $X_{i}\in E$ 
    vanishes whenever it contains $N$ times the same argument, that 
    is if there are $N$ distinct elements $i_{1},\dots, i_{N}$ of 
    $\{1,\dots,n\}$ such that $X_{i_{1}}=\dots = X_{i_{N}}$.

    \noindent {\bf Proof} This is straightforward, as for the proof 
    of Lemma 1, since one has more than $N$ symmetrized entries which 
    are distributed among less than $N-1$ columns in which the 
    entries are antisymmetrized.$\square$\\
 
    The right action $\lambda_{(Y^N)}$ of $\Tg(E)$ on $\wedge_{N}E$ will 
	   also be simply denoted by $\lambda_{N}$. In the case $N=2$, 
	   $\wedge_{2}E$ is the usual exterior algebra $\wedge E$ of $E$ and 
	   the right action $\lambda_{2}$ of $\Tg(E)$ factorizes through the 
	   right action of $\wedge E$ on itself, in particular 
	   $\ker(\lambda_{2})$ is the two-sided ideal of $\Tg(E)$ generated by 
	   the $X\otimes X$ for $X\in E$. Thus the graded algebra 
	   $\wedge_{[(Y)]}E=\Tg(E)/\ker \lambda_{(Y)}$ is also a 
	   generalization of the exterior algebra of $E$. For $(Y)=(Y^N)$, 
	   $\wedge_{[(Y^N)]}E$ will be simply denoted by $\wedge_{[N]}E$. One
	   clearly has $\wedge_{[2]}E=\wedge_{2}E=\wedge E$ for $N=2$.
	   In the case $N=3$, it can be shown 
	   that $\ker(\lambda_{3})$ is the two-sided ideal of $\Tg(E)$ 
	   generated by the
	   \[
	   X\otimes Y\otimes Z + Z\otimes X\otimes Y + Y\otimes Z \otimes X
	   \]
	   and the 
	   \[
	   X\otimes Y\otimes X\otimes X
	   \]
	   for $X,Y,Z\in E$. This implies that one has
	   \[
	   \lambda_{3}(X)\lambda_{3}(Y)\lambda_{3}(Z)+\lambda_{3}(Z)\lambda_{3}(X)\lambda_{3}(Y)+
	   \lambda_{3}(Y) \lambda_{3}(Z) \lambda_{3}(X)=0
	   \]
	   and
	   \[
	   \lambda_{3}(X)\lambda_{3}(Y)(\lambda_{3}(X))^2=0
	   \]
	   for any $X,Y,Z\in E$ and that these are the only independent 
	   relations in the associative algebra $\im 
	   (\lambda_{3})=\wedge_{[3]}E$. This means that $\wedge_{[3]}E$ is 
	   the associative unital algebra generated by the subspace $E$ with 
	   relations $XYZ+ZXY+YZX=0$ and $XYX^2=0$ for $X,Y,Z\in E$. The 
	   graduation is induced by giving the degree one to the elements of 
	   $E$ which is consistent since the relations are homogeneous. It 
	   is clear on this example that the homogeneous subspaces 
	   $\wedge^p_{[N]}E$ of $\wedge_{[N]}E$ are generally not 
	   irreducible for the (obvious) action of $GL(E)$. It is not hard 
	   to see that one has 
	       \[
	       \omega_{0}\wedge_{[N]}E=\wedge_{N}E
	       \]
	       where $\omega_{0}$ is a generator ($\simeq \bbbone$) of 
	       $\wedge^0_{N}E\simeq \mathbb R$, that is $\omega_{0}$ is 
	        a cyclic vector for the action of $\wedge_{[N]}E$ on 
	        $\wedge_{N}E$.\\
	       
	       Corresponding to the generalization $\wedge_{[(Y)]}E$ of the 
	       exterior algebra there is a generalization 
	       $\Omega_{[(Y)]}(M)$ of differential forms on a smooth 
	       manifold $M$ which is defined in a similar way as 
	       $\Omega_{(Y)}(M)$ was defined in Section 3. This 
	       $\Omega_{[(Y)]}(M)$ is then a graded associative algebra and 
	       $\Omega_{(Y)}(M)$ is a right graded 
	       $\Omega_{[(Y)]}(M)$-module (etc.). In the case $(Y)=(Y^N)$ 
	       one writes $\Omega_{[N]}(M)$ for this generalization. For 
	       $M=\mathbb R^D$ one has
	       \[
	       \Omega_{[N]}(\mathbb R^D)=\wedge_{[N]}\mathbb 
	       R^{D\ast}\otimes C^\infty (\mathbb R^D)
	       \]
	       and, by identifying $\Omega_{[N]}(\mathbb R^D)$ as a 
	       graded-subspace of $\Tg(\mathbb R^{D\ast})\otimes C^\infty 
	       (\mathbb R^D)$ and by using the canonical flat torsion-free 
	       linear connection of $\mathbb R^D$ one can define a 
	       $N$-differential $d$ on $\Omega_{[N]}(\mathbb R^D)$ by 
	       appropriate projection. One can proceed similarity for 
	       $\Omega_{[(Y)]} (\mathbb R^D)$ when $(Y)$ satisfies the 
	       assumption of Lemma 1 (or Lemma 2, Lemma 7). More precisely, the 
	       $N$-complexes constructed so far are particular cases of the 
	       following general construction.\\
	       
	       Let $\cala=\oplus_{n\in \mathbb N}\cala^n$ be an associative 
	       unital graded algebra generated by $D$ elements of degree 
	       one $\theta^\mu$ for $\mu\in \{ 1,\dots,D\}$ such that
	       \begin{equation}\label{Nil}
		\sum_{p\in \cals_{N}} \theta^{\mu_{p(1)}}\dots \theta^{\mu_{p(N)}}=0
		\end{equation}
		for any $\mu_{1},\dots,\mu_{N}\in \{ 1,\dots, D\}$. Then the 
		algebra $\cala(\mathbb R^D)$ defined by $\cala(\mathbb 
		R^D)=\cala\otimes C^\infty (\mathbb R^D)$ is a graded algebra and 
		one defines a $N$-differential $d$ on $\cala(\mathbb R^D)$, i.e. a 
		linear mapping $d$ of degree one satisfying $d^N=0$, by setting
		\begin{equation} \label{Nd}
		 d(a\otimes f)=(-1)^n a \theta^\mu \otimes \partial_{\mu} f
		 \end{equation}
		 for $a\in \cala^n$ and $f\in C^\infty(\mathbb R^D)$. Let 
		 $\calm=\oplus_{n}\calm^n$ be a 
		 graded right $\cala$-module, then $\calm(\mathbb R^D)=\calm \otimes 
		 C^\infty(\mathbb R^D)$ is a graded space which is a graded right 
		 $\cala(\mathbb R^D)$-module and one defines a $N$-differential $d$
		 on $\calm(\mathbb R^D)$ by setting
		 \begin{equation}\label{NdM}
		 d(m\otimes f)=(-1)^n m\theta^\mu \otimes \partial_{\mu}f
		 \end{equation}
		 for $m\in \calm^n$ and $f\in C^\infty(\mathbb R^D)$.  The 
		 (irrelevant) sign $(-1)^n$ in formulas (\ref{Nd}) and (\ref{NdM}) 
		 is here in order to recover the usual exterior differential in the 
		 case where $\cala=\wedge \mathbb R^{D\ast}=\calm$.\\
		 
		 It is clear that $\Omega_{[N]}(\mathbb R^D)=\cala(\mathbb R^D)$ 
		 for $\cala=\wedge_{[N]}\mathbb R^{D\ast}$ and that 
		 $\Omega_{N}(\mathbb R^D)=\calm(\mathbb R^D)$ for 
		 $\calm=\wedge_{N}\mathbb R^{D\ast}$. If $(Y)$ satisfies the 
		 assumption of Lemma 1 one can take (in view of Lemma 7) 
		 $\cala=\wedge_{[(Y)]}\mathbb R^{D\ast}$ and 
		 $\calm=\wedge_{(Y)}\mathbb R^{D\ast}$ and then 
		 $\Omega_{[(Y)]}(\mathbb R^D)=\cala(\mathbb R^D)$ and 
		 $\Omega_{(Y)}(\mathbb R^D)=\calm(\mathbb R^D)$.

\section{Further remarks}

Our original unpublished project for proving Theorem 1 was based on 
the construction of generalized algebraic homotopy in appropriate 
degrees. Let us explain what it means, why it is rather cumbersome and 
why the proof given here, based on the introduction of the multigraded 
differential algebra $\fraca$, is much instructive and general and is 
related to the ansatz of Green for the fermionic parastatistics of 
order $N-1$ (in the case $d^N=0$).\\

Let $\Omega=\oplus_{n} \Omega^n$ be a $N$-complex (of cochains) with 
$N$-differential $d$. An {\sl algebraic homotopy} for the degree $n$ 
will be a family of $N$ linear mappings
\[
h_{k} : \Omega^{n+k}\rightarrow \Omega^{n+k-N+1}
\]
for $k=0,\dots,N-1$ 
such that $\sum^{N-1}_{k=0}d^{N-1-k} h_{k}d^k$ is the identity 
mapping $I_{n}$ of $\Omega^n$ onto itself. If such a homotopy exists 
for the degree $n$, then one has $H^n_{(k)}=0$ for $k\in \{ 0,\dots, 
N-1\}$. Indeed let $\omega\in \Omega^n$ be such that $d^k\omega=0$ 
then one has $\omega=d^{N-k}\left (\sum^{k-1}_{p=0} d^{k-1-p} 
h_{p}d^p\omega \right)$.\\

Our original strategy for proving Theorem 1 was to show that one can 
construct inductively such homotopies for the degrees $(N-1)p$ with 
$p\geq 1$ in the case of the $N$-complex $\Omega_{N}(\mathbb R^D)$ and 
our idea was to exhibit explicit formulas. Unfortunately this latter task 
seems very difficult in general. We only succeeded in producing formulas 
in a closed form in the case $N=3$ and we refrain to give them here 
because this would imply explanations of our normalization 
conventions which have no character of naturality. The difficulty is 
indeed a problem of normalization. For the classical case $N=2$, one 
obtains a homotopy formula by using the inner derivation $i_{x}$ with 
respect to the vector field $x$ with components $x^\mu$. In this case 
one uses the fact that both $d$ and $i_{x}$ are antiderivations and 
that the Lie derivative $L_{x}=di_{x}+i_{x}d$ is the sum of the 
form-degree and of the degree of homogeneity in $x$. This gives 
homotopy formulas for forms which are homogeneous polynomials in $x$ 
and one gets rid of the above degree by appropriately weighted radial 
integration and obtains thereby the usual homotopy formula for 
positive form-degree. In this case the normalizations are fixed by 
the (anti)derivation properties. In the case $N\geq 3$, $d$ has no 
derivation property and one has to generalize $i_{x}$ which is 
possible with $i^N_{x}=0$ but there is no natural normalization since 
$i_{x}$ cannot possess derivation property. As a consequence the 
appropriate generalization of the Lie derivative involves a linear 
combination of products of powers of $d$ and $i_{x}$ with 
coefficients which have to be fixed at each step. That this is possible 
constitutes a cumbersome proof of Theorem 1 but does not allow easily 
to write closed formulas.\\

The interest of the proof of Theorem 1 presented here lies in the fact 
that it follows from the more general Theorem 2 which can be applied 
to other situations in particular to investigate the generalized 
cohomology of $\Omega_{[N]}(\mathbb R^D)$. Moreover, the realization 
of $\Omega_{N}(\mathbb R^D)$ embedded in $\fraca$ is related to the 
Green ansatz for the parafermionic creation operators of order $N-1$. 
Indeed if instead of equipping $\fraca$ with the graded commutative 
product one replaces in the definition of $\fraca$ the graded tensor 
products of graded algebras by the ordinary tensor products of 
algebras (applying the appropriate Klein transformation) then the 
$d_{i}x^\mu$ and the $d_{j}x^\nu$ commute for $i\not= j$ and the 
$d_{i}$ defined by the same formulas (14) commute, i.e. satisfy 
$d_{i}d_{j}=d_{j}d_{i}$ instead of (15), from which it follows that 
$\sum_{i}d_{i}$ is only a $N$-differential. This latter 
$N$-differential induces the $N$-differential of $\Omega_{N}(\mathbb 
R^D)\subset \fraca$ and the relation with the Green ansatz becomes 
obvious after Fourier transformation.\\

The basic $N$-complexes considered in this paper are $N$-complexes of 
smooth tensor fields on $\mathbb R^D$ and we have seen the difficulty 
to extend the formalism on an arbitrary $D$-dimensional manifold $M$. 
In the case of a complex (holomorphic) manifold $M$ of complex 
dimension $D$, there is an extension of the previous formalism at the 
$\bar \partial$-level which we now describe shortly.\\

Let $M$ be a complex manifold of complex dimension $D$ and let $T$ be 
a smooth covariant tensor field of type $(0,p)$ (i.e. of $d\bar 
z$-degree $p$) with local components $T_{\bar\mu_{1}\dots 
\bar\mu_{p}}$ in local holomorphic coordinates $z^1,\dots,z^D$. Then 
$\partial_{\bar\mu_{p+1}}T_{\bar\mu_{1}\dots 
\bar\mu_{p}}$  are the components of a well-defined smooth covariant 
tensor field $\bar\nabla T$ of type $(0,p+1)$ since the transition 
functions are holomorphic, where $\partial_{\bar\mu}$ denotes the 
partial derivative $\partial/\partial\bar z^\mu$ of smooth functions. 
Let $(Y)$ be a sequence $(Y_{p})_{p\in\mathbb N}$ of Young diagrams 
such that $\vert Y_{p}\vert=p$ ($\forall p\in \mathbb N$) and denote 
by $\Omega^{0,p}_{(Y)}(M)$ the space of smooth covariant tensor fields 
of type $(0,p)$ with Young symmetry type $Y_{p}$ (with obvious 
notation). Let us set 
$\Omega^{0,\ast}_{(Y)}(M)=\oplus_{p}\Omega^{0,p}_{(Y)}(M)$ and 
generalize the $\bar\partial$-operator by setting
\[
\bar\partial = (-1)^p \Yg_{p+1}\circ \bar \nabla : 
\Omega^{0,p}_{(Y)}(M)\rightarrow \Omega^{0,p+1}_{(Y)}(M)
\]
with obvious conventions. It is clear that if $(Y)$ is such that for 
any $p\in \mathbb N$ the number of columns of $Y_{p}$ is strictly less 
than $N$, then one has $\bar\partial^N=0$ so 
$\Omega^{0,\ast}_{(Y)}(M)$ is a $N$-complex (for $\bar\partial$). In 
particular one has the $N$-complex $\Omega^{0,\ast}_{N}(M)$ for 
$\bar\partial$ by taking $(Y)=(Y^N)$. One has an obvious extension of 
Theorem~1 ensuring that the generalized $\bar\partial$-cohomology of 
$\Omega^{0,\ast}_{N}(\mathbb C^D)$ vanishes in degree $(N-1)p$ (i.e. 
bidegree or type $(0,(N-1)p))$ for $p\geq 1$. It is thus natural to 
seek for an interpretation of this generalized cohomology for 
$\Omega^{0,\ast}_{N}(M)$ in degrees $(N-1)p$ with $p\geq 1$ for an 
arbitrary complex manifold $M$ and one may wonder whether it can be 
computed in terms of the ordinary $\bar\partial$-cohomology of $M$.

\end{document}